\input ./style/arxiv-general.cfg
\documentclass[aop,MSNbibl,seceqn,dvips]{arximspdf}
\makeatletter
   \@ifpackageloaded{graphicx}{}{\usepackage{graphicx}}
\makeatother
\doi{10.1214/15-AOP1036}% Updated by VTEXPTS2LaTeX.exe, 11.08.2015
\volume{44}
\issue{4}
\pubyear{2016}
\firstpage{2858}
\lastpage{2888}
\docsubty{FLA}

\makeatletter
\renewcommand{\det}{\operatorname{det}}
\newcommand{\rrvert}{\vert}
\newcommand{\llvert}{\vert}
\renewcommand{\mid}{|}
\newcommand{\rrVert}{\Vert}
\newcommand{\llVert}{\Vert}
\newtheorem{theorem}{Theorem}[section]
\newtheorem{proposition}[theorem]{Proposition}
\newtheorem{corollary}[theorem]{Corollary}
\newtheorem{lemma}[theorem]{Lemma}
\newproclaim{remark}[theorem]{Remark}
\newproclaim{definition}[theorem]{Definition}
\newproclaim{assumption}[theorem]{Assumption}
\makeatother

\begin{document}
\begin{frontmatter}

\title{Hafnians, perfect matchings and Gaussian matrices}
\runtitle{Hafnians, perfect matchings and Gaussian matrices}

\begin{aug}
% Corresponding author: Ofer Zeitouni - ofer.zeitouni@weizmann.ac.il% Updated by VTEXPTS2LaTeX.exe, 12.08.2015 15:12
%Updated by VTEXPTS2LaTeX.exe, 11.08.2015 13:35
\author[A]{\fnms{Mark}~\snm{Rudelson}\thanksref{M1,T1}\ead[label=e1]{rudelson@umich.edu}},
\author[B]{\fnms{Alex}~\snm{Samorodnitsky}\thanksref{M2,T2}\ead[label=e2]{salex@cs.huji.ac.il}}
\and
\author[C]{\fnms{Ofer}~\snm{Zeitouni}\corref{}\thanksref{M3,M4,T3}\ead[label=e3]{ofer.zeitouni@weizmann.ac.il}}
\runauthor{M. Rudelson, A. Samorodnitsky and O. Zeitouni}
\affiliation{University of Michigan\thanksmark{M1}, Hebrew University of Jerusalem\thanksmark{M2},
Weizmann Institute\thanksmark{M3} and New York University\thanksmark{M4}}
%\dedicated{}
\address[A]{M. Rudelson\\
Department of Mathematics\\
University of Michigan\\
530 Church Street\\
Ann Arbor, Michigan 48104\\
USA\\
\printead{e1}}
\address[B]{A. Samorodnitsky\\
School of Engineering and Computer Science\\
Hebrew University of Jerusalem\\
Givat Ram, Jerusalem 9190401\\
Israel\\
\printead{e2}}
\address[C]{O. Zeitouni\\
Faculty of Mathematics\\
Weizmann Institute\\
POB 26, Rehovot 76100\\
Israel\\
and\\
Courant Institute\\
New York Universirty\\
251 Mercer Street\\
New York, New York 10012\\
USA\\
\printead{e3}}
\end{aug}
\thankstext{T1}{Supported in part by NSF Grant DMS-11-61372 and USAF
Grant FA9550-14-1-0009.}
\thankstext{T2}{Supported in part by grants from the US--Israel
Binational Science Foundation and from the Israel Science Foundation.}
\thankstext{T3}{Supported in part by a grant from the Israel Science
Foundation and by the
Herman P. Taubman professorial chair of Mathematics at WIS.}

% HISTORY:
%
\received{\smonth{9} \syear{2014}}% Updated by VTEXPTS2LaTeX.exe,
%11.08.2015 13:35
%
\revised{\smonth{6} \syear{2015}}% Updated by VTEXPTS2LaTeX.exe,
%11.08.2015 13:35

% ABSTRACT
%
\begin{abstract}
We analyze the behavior of the Barvinok estimator of the hafnian of
even dimension, symmetric
matrices with nonnegative entries. We introduce a condition under
which the Barvinok estimator achieves subexponential errors, and show
that this condition is almost optimal. Using that hafnians count the
number of perfect matchings in graphs, we conclude that Barvinok's
estimator gives a polynomial-time algorithm for the approximate (up to
subexponential errors) evaluation of the number of perfect matchings.
\end{abstract}

% KEYWORDS
\end{frontmatter}

%s1 #&#
\section{Introduction}\label{sec1}
The number of perfect matchings in a bipartite graph is given by the
permanent of the bipartite adjacency matrix of the graph. Since
computing the
permanent is generally computationally hard \cite{Val}, various algorithms
have been proposed to compute it approximately. We mention in
particular the
MCMC algorithm of Jerrum--Sinclair--Vigoda \cite{JSV}, the
Linial--Samorodnitsky--Wigderson rescaling algorithm \cite{LSW}
(denoted LSW in the sequel), and the Barvinok--Godsil--Gutman algorithm
\cite{GG,Bar}; the analysis of the latter algorithm was the subject of
%Ofer #1
the previous work \cite{RZ}.

{A more general (and hence hard)} combinatorial problem is that of
computing the number of perfect matchings in a graph with an even
number of vertices. Let $A$ denote the
adjacency matrix of such a graph with $n=2m$ vertices.
The relevant combinatorial notion here is the \textit{hafnian} \cite
{Minc}, defined as
\[
\operatorname{haf}(A)=\frac{1}{m! 2^m}\sum_{\sigma\in\mathcal
{S}_{n}}
\prod_{j=1}^m A_{\sigma(2j-1),\sigma(2j)},
\]
where $\mathcal{S}_{n}$ denotes the symmetric group on $[n]$. It is
immediate to
check (see, e.g.,~\cite{Bar}) that
%
%e1.1 #&#
%
\begin{equation}
\label{haffnian-matching} \#\mbox{perfect matchings in $A$} = \operatorname{haf}(A).
\end{equation}
Thus, the interest in an efficient computation of $\operatorname{haf}(A)$.
As for the permanent, the exact computation of $\operatorname{haf}(A)$ is
computationally expensive.
This problem of estimating the hafnian
%\textcolor{blue}
{seems to be harder} to attack than the corresponding problem for the
permanent since many algorithms known for permanent
%Ofer #2
approximation
break down when extended to hafnians.
In particular, the LSW rescaling algorithm \cite{LSW} transforms the
adjacency matrix of a graph to an almost doubly stochastic one. Yet, a
nontrivial lower estimate of the hafnian of a doubly stochastic matrix
is impossible; see \cite{BS}.
Also, in contrast with the computation of the permanent,
\cite{JSV} points out that the proof of convergence of the MCMC
algorithm breaks down for
the approximate computation of the hafnian (unless the minimal degree
is at least $n/2$, see \cite{JS}).

We consider in this paper
the computation of $\operatorname{haf}(A)$ for symmetric
matrices with nonnegative
entries. Note that
the diagonal entries play no role in the computation of
$\operatorname{haf}(A)$ and, therefore, in the rest of this paper we
always assume
that $A_{ii}=0$ for all $i$.

In his seminal paper \cite{Bar}
discussing the Godsil--Gutman estimator for the permanent,
Barvinok also introduces a probabilistic estimator of $\operatorname
{haf}(A)$ for a
symmetric matrix $A$
possessing
nonnegative entries. Let $W$ be a real skew symmetric matrix with
independent centered normal entries $W_{ij}$ above the diagonal
satisfying $\mathbb{E}W_{ij}^2=A_{ij}$.
In\vspace*{2pt} other words, let
$G^{\mathrm{skew}}$ denote a skew symmetric matrix with independent $N(0,1)$
entries above the main diagonal.
Let $W=W(A)$ denote the skew symmetric
matrix
with
%
%e1.2 #&#
%
\begin{equation}
\label{Wdef} W_{ij}=G_{ij} \sqrt{A_{ij}},
\qquad i<j
\end{equation}
and write $W={\mathcal{A}}\odot G^{\mathrm{skew}}$,
where
${\mathcal{A}}$ denotes the element-wise square-root of $A$, that is,
${\mathcal{A}}_{ij}=\sqrt{A_{ij}}$.
Then
%
%e1.3 #&#
%
\begin{equation}
\label{eq-determinantrelation} \operatorname{haf}(A)=E\det(W).
\end{equation}
Thus, $\det(W)$, which is an easily computable quantity, is a consistent
estimator for $\operatorname{haf}(A)$, and Barvinok \cite{Bar}
proceeds to prove
that for any matrix $A$, $e^{-\gamma n} \operatorname{haf}(A) \leq
\det(W)\leq
C\cdot\operatorname{haf}(A)$ with high probability, where $\gamma$
is Euler's constant.
Other approaches to computing the hafnian include \cite{BS} (which
however does not apply to adjacency matrices of nontrivial graphs),
\cite{BGKNT}, where a deterministic algorithm of %\textcolor{blue}
{subexponential} complexity
is constructed and analyzed, and \cite{BS1}, where a random algorithm
is analyzed but the precision of the algorithm depends in a complicated
way on the number of perfect matchings.
%\marginpar{check references}

%Our goal in this paper is to analyze the performance of the Barvinok
%estimator
%for the hafnian. On a technical level, the main difference between this
%analysis and the one performed in \cite{RZ} is that the Gaussian matrix
%$W$ is skew-symmetric, and thus with less independence than the
%unrestricted
%matrix $W$ in \cite{RZ}. Consequently, the method for estimating the
%singular values of $W$ differs from that in \cite{RZ}: we modify the
%argument developed there for estimating the minimal singular value so
%that it
%can handle the extra dependencies due to the skew-symmetry. On the
%other hand, instead of developing a estimate for intermediate singular
%values
%as in \cite{RZ} (a difficult task here due to the skew-symmetry), we
%use
%the fact that $iW$ is Hermitian, which allows us
%to use estimates from recent work \cite{EKYY} on the local semi-circle
%law.

%MR
Our goal in this paper is to analyze the performance of the Barvinok estimator
for the hafnian.
As in \cite{RZ}, establishing the concentration of a random
determinant hinges on bounding the singular values of the Gaussian
matrix $W$. This crucial step, however, essentially differs from \cite
{RZ} as $W$ is skew-symmetric, and thus has less independence than the
unrestricted matrix in \cite{RZ}. Handling these dependences required
different arguments for the smallest and the intermediate singular
values. In the first case, we employ a conditioning argument tailored
to take into account the structure of the graph (Lemmas~\ref{l: tree},
\ref{l: moderately compressible individual}). The fact that the
entries of $W$ are real, and thus the $(n-1) \times(n-1)$ main minors
of it are degenerate plays a central role here. On the
other hand, instead of developing a estimate for intermediate singular values
as in \cite{RZ} (a difficult task here due to the skew-symmetry), we use
the fact that the imaginary-valued matrix $iW$ is Hermitian, which
allows us
to use estimates from recent work \cite{EKYY} on the local semi-circle law.
% MR end

To formulate our results, we introduce a notion of strong expansion for
graphs. This notion strengthens the standard notion of vertex expansion
assuming that sets having many connected components expand faster. For
a set $J \subset[n]$ of vertices, denote by $\operatorname{Con}(J)$
the set of the connected
components of $J$,
and by $\partial(J)$ the boundary of $J$, that is,
$\partial(J)=\{i\in[n]\setminus J: \exists j\in J, (i,j)$ is an edge$\}$.
%
%de1.1 #&#

\begin{definition}
\label{def-1}
Let $\kappa \in(0,1)$, and let $1<m<n$.
We say that the graph $\Gamma $ is strongly expanding with parameter
$\kappa $
up to level $m$ if for any set $J \subset[n]$ of vertices with $
\llvert  J\rrvert
\le m$,
\[
\bigl\llvert \partial(J) \bigr\rrvert - \bigl\llvert \operatorname {Con}(J)
\bigr\rrvert \ge\kappa \cdot\llvert J \rrvert .
\]
%
% We will abbreviate this property by writing $\G\in\mbox{Ex}(\k,m)$.
\end{definition}

In this definition and below, we use the following notational
convention. Important parameters which appear in definitions and
theorems are denoted by Greek letters. Unimportant constants whose
value may change from line to line are denoted $c,c', C$, etc.

The simplest form of our results is in case $A$ is the adjacency matrix
of a $d$-re\-gular graph.
%
%th1.2 #&#

\begin{theorem}
\label{thm-3}
Fix $\alpha ,\kappa>0$. Let $A$ be the adjacency matrix of a $d$-regular
graph $\Gamma$ with $d\geq\alpha n+2$. Assume that
%
%e1.4 #&#
%
\begin{equation}
\label{eq-SEthm12} \mbox{$\Gamma$ is $\kappa$ strongly expanding up to level $n(1-
\alpha )/(1+\kappa/4)$.}
\end{equation}
%
%\begin{longlist}[(2)]
%

(1) Then for any $\epsilon<1/5$ and $D>4$,
%
%e1.5 #&#
%
\begin{equation}
\label{eq-thm3-1bis} \mathbb{P} \bigl( \bigl\llvert \log\operatorname {haf}(A)-\log\det
(W) \bigr\rrvert > C n^{1-\epsilon} \bigr)\leq n^{-D},
\end{equation}
where $C=C(\alpha ,\kappa,D,\epsilon)>0$.

(2) Fix $\delta>0$.
If, in addition to the above assumptions, the matrix $A/d$
possesses a spectral gap $\delta$, then for any $D>4$,
%Ofer #3
%there exists
%$C=C(\a,\kappa,D,\delta)>0$ so that
%
%e1.6 #&#
%
\begin{equation}
\label{eq-thm3-1} \mathbb{P} \bigl( \bigl\llvert \log\operatorname {haf}(A)-\log\det
(W) \bigr\rrvert >C n^{1/2}\log^{1/2} n \bigr)\leq
n^{-D},
\end{equation}
where
$C=C(\alpha ,\kappa,D,\delta)>0$.
%\end{longlist}
%
\end{theorem}

[The assumption in (2) means that the modulus of the
eigenvalues of $A/d$ is either $1$ or smaller than $1-\delta$.]
Theorem~\ref{thm-3} is an immediate consequence of our more general
Theorem~\ref{thm-2} below.

We discuss the definition of strongly expanding graphs in
Remark~\ref{r: example} below.

The extension of Theorem~\ref{thm-3}
to irregular graphs requires the notion of doubly stochastic scaling of
matrices.
We also need the notion of spectral gap for stochastic matrices.
%
%de1.3 #&#

\begin{definition}
\label{def-ds}
A matrix $A$ with nonnegative entries is said to possess a doubly
stochastic scaling if there exist two diagonal matrices $D_1,D_2$ with
positive entries such that the matrix $B=D_1AD_2$ is doubly stochastic,
that is, $\sum_j B_{ij}=\sum_k B_{k\ell}=1$ for all $i,\ell$. We
call such $B$ a doubly stochastic scaling of $A$.
\end{definition}

%de1.4 #&#

\begin{definition}
\label{def-spectralgap}
A symmetric stochastic matrix $A$ is said to possess
a \textit{spectral gap} $\delta$ if
there do not exist eigenvalues of $A$ in $(-1,-1+\delta)\cup(1-\delta,1)$.
\end{definition}

We will show below (see Corollary~\ref{cor: matching exists}) that
the adjacency matrix of a strongly expanding graph with appropriate
lower bound on its minimal degree
possesses a unique doubly stochastic scaling, with $D_1=D_2$.
We use this fact in the following theorem, where the
required lower
bound on the minimal degree is satisfied.
%
%th1.5 #&#

\begin{theorem}
\label{thm-3bis}
Fix $\alpha ,\kappa,\vartheta>0$. Let $A$ be the adjacency matrix of
a graph
$\Gamma$ whose minimal degree satisfies $d\geq\alpha n+2$.
%\textcolor{red}{Let $B$ denote the doubly stochastic scaling of $A$.}
Assume that
%
%e1.7 #&#
%
\begin{equation}
\label{eq-SEthm12bis} \mbox{$\Gamma$ is $\kappa$ strongly expanding up to level $n(1-
\alpha )/(1+\kappa/4)$,}
\end{equation}
and
%
%e1.8 #&#
%
\begin{equation}
\label{eq-dsadj} \mbox{the doubly stochastic scaling $B$ of $A$ satisfies $\max
_{i,j} B_{ij}\leq n^{-\vartheta}$}.
\end{equation}
%
%\begin{longlist}[(2)]

(1) Then for any $\epsilon<1/5$ and $D>4$,
%
%e1.9 #&#
%
\begin{equation}
\label{eq-thm3-1newbis} \mathbb{P} \bigl( \bigl\llvert \log \operatorname{haf}(A)-\log
\det (W) \bigr\rrvert > C n^{1-\epsilon\vartheta} \bigr)\leq n^{-D}
\end{equation}
with $C=C(\alpha ,\kappa,D,\epsilon)>0$.

(2) Fix $\delta>0$. If, in addition to the above assumptions,
$B$ possesses a spectral gap $\delta$,
then
for any $D>4$,
%Ofer #4 there exists
%$C=C(\a,\kappa,D,\delta, \vartheta)>0$ so that
%
%e1.10 #&#
%
\begin{equation}
\label{eq-thm3-1new} \mathbb{P} \bigl( \bigl\llvert \log\operatorname {haf}(A)-\log
\det (W) \bigr\rrvert >C n^{1-\vartheta/2}\log^{1/2} n \bigr)\leq
n^{-D},
\end{equation}
where
$C=C(\alpha ,\kappa,D,\delta, \vartheta)>0$.
%\end{longlist}
%
\end{theorem}

Condition\vspace*{1pt} (\ref{eq-dsadj}) can be readily checked in polynomial time
by applying the LSW scaling algorithm, stopped when its error is
bounded above by $n^{-1}$. Indeed, at such time, the LSW algorithm
output is a matrix $C=C(A)$ which is almost doubly stochastic in the
sense that, with $B=B(A)$ denoting the doubly stochastic scaling of
$A$, one has $\max_{ij}\llvert    B(A)-C(A)\rrvert   <n^{-1}$. Because the
maximal entry
of $B$ is at least $n^{-1}$, this implies that the maximal entries of
$B$ and of $C$ are of the same order.
Note also that the spectral gap condition in point (2) of
Theorem
\ref{thm-3bis}, which depends only on the eigenvalues of $B$, can also be
checked in polynomial time.

%\textcolor{blue}
{We note that for a given $0 < \vartheta< 1$, there exist stronger
expansion conditions on the graph $\Gamma$ which ensure that the
maximal element in the doubly stochastic scaling of its adjacency
matrix is of size at most $n^{-\vartheta}$. That is, if $\Gamma$
satisfies these stronger properties, condition (\ref{eq-dsadj}) is
automatically satisfied. We refer to Section~\ref{sec: doubly
stochastic}, Proposition~\ref{prop:not-large}, for details.}

Conditions
%\eqref{eq-SEpropds}
(\ref{eq-SEthm12bis}) and (\ref{eq-dsadj}) play different roles in
the proof. The first one is needed to establish the lower bound on the
smallest singular value of $W$, and the second one guarantees that most
of the singular values are greater than $n^{-\epsilon}$.

%re1.6 #&#
%
\begin{remark} \label{r: example}
The definition of strongly expanding graphs (Definition~\ref{def-1}
above) reminds one of that of a vertex expander. Yet, it is stronger in
two senses. First, the strong expansion property takes into account the
geometry of the set, requiring more rapid expansion for more ``spread
out'' sets. Second, we want this expansion property to hold for all
sets of size relatively close to $n$, while for the classical
expanders, the corresponding property is required only for sets with at
most $n/2$ vertices. This may look unnatural at first glance. However,
% \textcolor{red}{Modify after counter example has stabilized: However,
%Proposition %\ref{p: graph-example} below shows that
one may construct an example of a graph that can have the strong
expansion property up to a level arbitrary close to $1$, and yet the
matrix corresponding to it may be degenerate with probability~1.
Further, in Proposition
\ref{p: bias} below, we construct a graph whose adjacency matrix
barely misses the condition in Definition~\ref{def-1} and yet $\det
(W)/\operatorname{haf}(A)\leq e^{-cn}$ with high probability for an
appropriate $c>0$.
\end{remark}

%pr1.7 #&#
%
\begin{proposition} \label{p: bias}
Let $\delta>0$. For any $N \in\mathbb{N}$, there exists a graph
$\Gamma $ with
$M>N$ vertices such that
%
%e1.11 #&#
%
\begin{equation}
\label{c: weak expansion} \qquad\forall J \subset[M] \qquad\llvert J\rrvert \le M/2\quad
\Rightarrow \quad \bigl\llvert \partial(J) \bigr\rrvert - (1-\delta) \bigl\llvert
\operatorname{Con}(J) \bigr\rrvert \ge\kappa \llvert J\rrvert
\end{equation}
and
\[
\mathbb{P} \biggl(\frac{\det(W)}{\mathbb{E}\det(W)} \le e^{-cM} \biggr) \ge1-
e^{-c'M}.
\]
Here, $c,c', \kappa $ are constants depending on $\delta$.
\end{proposition}

Our theorems on adjacency matrices are based on a general result
pertaining to doubly stochastic symmetric matrices $B$ with
nonnegative entries. We will consider matrices which have many
relatively large entries. To formulate this requirement precisely, we
introduce the notion of large entries graph.
%
%de1.8 #&#

\begin{definition} \label{def: large entries}
Let $A$ be a symmetric matrix with nonnegative entries.
For a parameter $\theta>0$, define the \emph{large entries graph}
$\Gamma
_A(\theta)$ by connecting the vertices $i,j \in[n]$ whenever $A_{ij}>
\theta$. If $A$ is the matrix of variances of entries of a skew
symmetric matrix $W$, we will also refer to $\Gamma _A(\theta)$ as the
\emph{large variances graph} of~$W$.
%\textcolor{red}{(or of $A$)}.
\end{definition}

We will now formulate two theorems on the concentration of the hafnian
of a skew symmetric matrix whose large variances graph satisfies a
strong expansion condition.
%
%th1.9 #&#

\begin{theorem}
\label{thm-1}
Fix $\beta ,\alpha ,\vartheta,\kappa>0$. Let $B$ be a symmetric
stochastic matrix of even size $n$ with nonnegative entries,
let $W=W(B)$ be as (\ref{Wdef}) (with $B$ replacing~$A$)
and
let $\Gamma =\Gamma _{B}(n^{-\beta })$ denote its large variances graph.
Assume that:
\begin{longlist}[(2)]
\item[(1)] The minimal degree of a vertex of $\Gamma $ is at least $\alpha n+2$.
\item[(2)] $\Gamma $ is $\kappa$ strongly-expanding up to level
$n(1-\alpha
)/(1+\kappa/4)$.
\item[(3)] $\max B_{ij}\leq n^{-\vartheta}$.
%\item$\Gamma$ possesses a perfect matching, that is, there exists a
%permutation %$\sigma\in\SS_{2n}$ so that $\prod_{i=1}^n \Gamma_{
%\sigma(2i-1),\sigma(2i)}=1$.
%\mbox{\rm\textcolor{red}{Check whether strong expansion implies}}
%\mbox{\rm\textcolor{red}{existence of perfect matching}}
\end{longlist}
%
% Let $W$ be an $n \times n$ skew symmetric random matrix with
%independent %up to the symmetry restriction centered normal entries
%satisfying $\E W_{ij}^2= %A_{ij}$.
Then, for any $\epsilon<1/5$ and $D>4$ there exists $C=C(\beta
,\alpha
,\kappa
,\vartheta,\epsilon,D)>0$ so that
%
%e1.12 #&#
%
\begin{equation}
\label{eq-thm1-1} \mathbb{P} \bigl( \bigl\llvert \log\operatorname {haf}(B)-\log
\det (W) \bigr\rrvert >C n^{1-\epsilon\vartheta} \bigr)\leq n^{-D}.
\end{equation}
%
%and
%\begin{equation}
%\label{eq-thm1-2}
%\P(\left\left\vert  \log\haf(A)-\log\det(W^2)\right\right\vert  >tn^{4/5})\leq
%e^{-c'' t^{1/5}}+
%\exp(-cn^{1/5}).
%\end{equation}
\end{theorem}

Somewhat tighter bounds are available if the matrix $B$
possesses a spectral gap.
%
%th1.10 #&#

\begin{theorem}
\label{thm-2}
Assume the conditions of Theorem~\ref{thm-1} and in addition assume that
the matrix $B$ has a spectral gap $\delta$. Then,
for any $D>4$,
%
%e1.13 #&#
%
\begin{equation}
\label{eq-thm2-1} \mathbb{P} \bigl( \bigl\llvert \log\operatorname {haf}(B)-\log
\det (W) \bigr\rrvert >C n^{1-\vartheta/2}\log^{1/2} n \bigr)\leq
n^{-D}.
\end{equation}
The constant $C$ here depends on all relevant parameters $\beta ,
\alpha , \kappa ,
\vartheta, \delta$ and $D$.
\end{theorem}

The structure of the paper is as follows. In Section~\ref{sec:
compressible}, we consider unit vectors that are close to vectors with
small support and derive uniform small ball probability estimates for
their images under the action of $W$. These estimates are used in
Section~\ref{sec: smallest singular} to obtain a lower bound for the
smallest singular values of $W$. In Section~\ref{sec: local bounds},
we provide local estimates for the empirical measure of eigenvalues of
$W$. Section~\ref{sec: proof thm 1 2} is devoted to the proof of
Theorems~\ref{thm-1} and~\ref{thm-2}. Section~\ref{sec: doubly
stochastic} is devoted to the proof of a combinatorial lemma concerning
the doubly stochastic scaling of adjacency matrices of strongly
expanding graphs,
which then is used in the proof of Theorem~\ref{thm-3bis}; in the
section we also present sufficient conditions that ensure that~(\ref
{eq-dsadj}) holds. Finally, in Section~\ref{sec: strong expansion} we
present the construction of the graph discussed in Proposition~\ref{p:
bias}, and provide the proof of the latter.

\section{Compressible vectors}
\label{sec: compressible}
To establish the concentration for the determinant of the matrix $W$,
we have to bound its smallest singular value. As is usual in this
context, we view the smallest singular value of a matrix as the minimum
of the norms of the images of unit vectors:
\[
s_n(W)=\min_{x \in S^{n-1}} \llVert Wx \rrVert
_2.
\]
Before bounding the minimal norm over the whole sphere, let us consider
the behavior of $\llVert  Wx \rrVert _2$ for a fixed $x \in S^{n-1}$.
We begin with a small ball probability estimate, which is valid for any
unit vector.
%
%le2.1 #&#

\begin{lemma} \label{l: sparse vector}
Let $W$ be a skew-symmetric $n \times n$ matrix with independent, up to
the symmetry restriction, normal entries.
Assume that for any $j \in[n]$, there exist at least $d$ numbers $i
\in n$ such that $\operatorname{Var}(w_{ij}) \ge n^{-c}$.
Then
for any $x \in S^{m-1}$, and for any $t>0$
\[
\mathbb{P} \bigl(\llVert Wx \rrVert _2 \le t n^{-c'} \bigr)
\le(Ct)^{- d},
\]
where $C,c'$ depend on $c$ only.
\end{lemma}

\begin{pf}
Let $x \in S^{n-1}$. Choose a coordinate $j \in[n]$ such that $|
x_j|
\ge n^{-1/2}$ and set $I=\{i \in[n] |  \operatorname{Var}(w_{ij})
\ge n^{-c} \}$. Condition on all entries of\vadjust{\goodbreak} the matrix $W$, except
those in the $j$th row and column. After this conditioning, for any $i
\in I$ the $i$th coordinate of the vector $Wx$ is a normal random
variable with variance $\operatorname{Var}(w_{i,j}) x_j^2 \ge
n^{-2c-1}$. Since the coordinates of this vector are conditionally
independent, an elementary estimate of the Gaussian density yields for
any $t>0$
\[
\mathbb{P} \bigl(\llVert Wx \rrVert _2< t n^{-c-1/2} |
w_{ab}, a,b \in[n] \setminus\{ j\} \bigr) \le(Ct)^{\llvert    I\rrvert   }.
\]
By the assumption of the lemma, $\llvert    I\rrvert   \ge d$.
Integration with respect to the other variables completes the proof.
\end{pf}

The next lemma is a rough estimate of the norm of a random matrix.
%
%le2.2 #&#

\begin{lemma} \label{l: norm}
Let $W$ be a be a skew-symmetric $n \times n$ matrix with independent,
up to the symmetry restriction, normal entries.
Assume that for any $i,j \in[n]$, $\operatorname{Var}(w_{ij}) \le1$.
Then
\[
\mathbb{P} \bigl(\llVert W \rrVert \ge n \bigr) \le e^{-n}.
\]
\end{lemma}

Lemma~\ref{l: norm} follows from the estimate $\llVert  W \rrVert ^2 \le\llVert  W \rrVert _{\mathrm{HS}}^2=\sum_{ij=1}^n
W_{ij}^2$, where the right-hand side is the
sum of squares of independent centered normal variables whose variances
are uniformly bounded.

Of course, the estimate in Lemma~\ref{l: norm} is very rough, but we
can disregard a constant power of $n$ in this argument. Lemma~\ref{l:
norm} allows us to extend the lower bound on the small ball probability
from a single vector to a neighborhood of a small-dimensional subspace.
To formulate it precisely, recall the definition of compressible and
incompressible vectors from \cite{RV,RV2}.
%
%de2.3 #&#

\begin{definition}
For $m <n$ and $v<1$, denote
%MathOfer
%
\begin{eqnarray*}
\operatorname{Sparse}(m) &=& \bigl\{ x \in S^{n-1} : \bigl\llvert
\operatorname{supp}(x) \bigr\rrvert \le m \bigr\}
\end{eqnarray*}
and
\begin{eqnarray*}
\operatorname{Comp}(m,v) &=& \bigl\{ x \in S^{n-1} :\exists y \in
\operatorname{Sparse}(m), \llVert x-y \rrVert _2 \le v \bigr\};
\\
\operatorname{Incomp}(m,v) &=& S^{n-1} \setminus\operatorname{Comp}(m,v).
\end{eqnarray*}
\end{definition}

The next lemma uses a standard net argument to derive the uniform
estimate for highly compressible vectors.
%
%le2.4 #&#

\begin{lemma} \label{l: very sparse}
Let $A$ be an $n \times n$ matrix satisfying the conditions of Lemmas
\ref{l: sparse vector} and~\ref{l: norm}.
Then
\[
\mathbb{P} \bigl( \exists x \in\operatorname{Comp} \bigl(d/2, n^{-c}
\bigr): \llVert Wx \rrVert _2 \le n^{-\bar{c}}\mbox{ and } \llVert
W \rrVert \le n \bigr) \le e^{- d/2},
\]
where $\bar{c}$ depends on $c$ only.
\end{lemma}

\begin{pf}
Let $t>0$ be a number to be chosen later, and set
\[
\varepsilon= t n^{-c'-2},
\]
where $c'$ is the constant from Lemma~\ref{l: sparse vector}.
Then there exists an $\varepsilon$-net $\mathcal{N}\subset
\operatorname
{Sparse}(d/2)$ of cardinality
\[
\llvert \mathcal{N}\rrvert \le\pmatrix{n
\cr
d/2} \cdot(3/\varepsilon)^{d/2}
\le \biggl( \frac{C n }{ t n^{-c'-2}} \biggr)^{d/2}.
\]
By Lemma~\ref{l: sparse vector} and the union bound,
\[
\mathbb{P} \bigl( \exists y \in\mathcal{N}:\llVert Wy \rrVert _2 \le
t n^{-c'} \bigr) \le \biggl( \frac{C n}{ t n^{-c'-2}} \biggr)^{d /2}
\cdot(Ct)^{d} \le e^{- d/2},
\]
provided that $t = n^{-c''}$ for an appropriately chosen $c''>0$.

Assume that for any $y \in\mathcal{N}, \llVert  Wy \rrVert _2
\ge t
n^{-c'}$. Let $x
\in\operatorname{Comp}(\delta/2, \varepsilon)=\break \operatorname
{Comp}(\delta/2, n^{-c})$, and choose $y \in\mathcal{N}$ be such that
$\llVert  x-y \rrVert _2< 2 \varepsilon$. If $\llVert  W
\rrVert \le n$, then
\[
\llVert Wx \rrVert _2 \ge\llVert Wy \rrVert _2- \llVert
W \rrVert \cdot \llVert x-y \rrVert _2 \ge t n^{-c'}- n
\cdot2 t n^{-c'-2} \ge n^{-\bar{c}}.
\]\upqed
\end{pf}

Our next goal is to show that the small ball probability estimate
propagates from strongly compressible vectors to moderately
compressible ones. At this step, the assumption that the large
variances graph is strongly expanding plays a crucial role. The strong
expansion condition guarantees that the matrix $W$ has enough
independent entries to derive the small ball estimate for a single
vector, despite the dependencies introduced by the skew-symmetric
structure. The next simple lemma is instrumental in exploiting the independence
that is still present.
%
%le2.5 #&#

\begin{lemma} \label{l: tree}
Let $T=(V,E)$ be a finite tree with the root $r \in V$. Assume that to
any $e \in E$ there corresponds a random variable $X_e$, and these
variables are independent. Assume also that to any $v \in V$ there
corresponds an event $\Omega _v$, which depends only on those $X_e$ for
which $v \in e$. Suppose that for any $v \in V$ and any $e_0$ connected
to $v$,
\[
\mathbb{P} \bigl(\Omega _v \mid \{X_e
\}_{e\neq e_0} %\mbox{all } e \neq e_0
 \bigr) \le p_v
\]
for some numbers $p_v \le1$. Then
\[
\mathbb{P} \biggl( \bigcap_{v \in V \setminus\{r\}} \Omega
_v \biggr) \le\prod_{v \in V \setminus\{r\}}
p_v.
\]
\end{lemma}

\begin{pf}
We prove this lemma by induction on the depth of the tree. Assume first
that the tree has depth 2. Then the statement of the lemma follows from
the fact that the events $\Omega _v, v \in V \setminus\{r\}$ are independent.

Assume now that the statement holds for all trees of depths smaller
than $k>2$ and let $T$ be a tree of depth $k$. Let $V_r$ and $E_r$ be
the sets of all vertices and edges connected to the root of the tree.
Then the events $\Omega _v, v \in V \setminus\{r\}$ conditioned on $x_e,
e \notin E_r$ are independent. Therefore,
\[
\mathbb{P} \biggl( \bigcap_{v \in V \setminus\{r\}} \Omega
_v \biggr) = \mathbb{E}\mathbb{P} \biggl[ \bigcap
_{v \in V \setminus\{
r\}} \Omega _v \Big| %e \notin E_r \right]
\{X_e \}_{e\notin E_r} \biggr] \le\prod
_{v \in V_r} p_v \cdot\mathbb{P} \biggl( \bigcap
_{v \in V \setminus
(V_r \cup\{r\})} \Omega _v \biggr).
\]
Note that the vertices $v \in V \setminus\{r\}$ form a forest with
roots $v \in V_r$. Since the events $\bigcap_{v \in T_l} \Omega _v$ are
independent for different trees $T_l$ in the forest, the statement of
the lemma follows by applying the induction hypothesis to each tree.
\end{pf}

Using Lemma~\ref{l: tree} and the strong expansion property of the
large variances graph, we establish the small ball probability bound
for the image of an incompressible vector.

%le2.6 #&#

\begin{lemma} \label{l: moderately compressible individual}
Let $c>0$.
Let $W$ be an $n \times n$ skew-symmetric centered Gaussian matrix.
Assume that its large variances graph $\Gamma _W(n^{-c})$ satisfies the
strong expansion condition with parameter $\kappa >0$ up to level $m<n$.
Let $t>0, k \le m$, and $v>0$.
Then for any $x \in\operatorname{Incomp}(k, v)$,
\[
\mathbb{P} \bigl(\llVert Wx \rrVert _2 \le n^{-C}v \cdot t
\bigr) \le t^{(1+\kappa ) k},
\]
where $C$ depends on $c$ only.
\end{lemma}

\begin{pf}
For $i \in[n]$, define the event $\Omega _i$ by
\[
\Omega _i= \bigl\{W: \bigl\llvert (W x)_i \bigr
\rrvert \le t n^{-(c+1)/2}v t \bigr\}.
\]
Let $J(x)=\{j: \llvert    x_j\rrvert   \ge n^{-1/2}v\}$. Since $x \in
\operatorname
{Incomp}(k, v)$, $\llvert    J(x)\rrvert   \ge k$.
Indeed, let $y \in\mathbb{R}^n$ be the vector containing $k$ largest in
absolute value coordinates of $x$. If $\llvert    J(x)\rrvert   \le k$, then
\[
\operatorname{dist} \bigl(x, \operatorname{Sparse}(k) \bigr) \le \llVert x-y
\rrVert _2 \le \biggl( \sum_{j \notin J(x)}
x_j^2 \biggr)^{1/2} \le v.
\]

Choose a subset $J \subset J(x)$ with $\llvert  J\rrvert  =k$. For $i
\in[n]$ set
$p_i=t$ whenever $i \sim j$ for some $j \in J$;
otherwise set $p_i=1$. Then for any $j_0 \in J$ and for any $i_0 \sim j_0$,
\[
\mathbb{P} \bigl(\Omega _{i_0} \mid W_{ij}, (i,j)
\neq(i_0,j_0) \bigr) \le t.
\]
Indeed, $(Wx)_i$ is a normal random variable with variance at least
\[
\operatorname{Var}(w_{i_0 j_0}) \cdot x_{j_0}^2 \ge
n^{-c} \cdot n^{-1} v^2,
\]
so the previous inequality follows from the bound on the maximal density.

To prove Lemma~\ref{l: moderately compressible individual},
we will use Lemma~\ref{l: tree}. To this end, we will construct a
forest consisting of $L=\llvert   \operatorname{Con}(J)\rrvert   $ trees with
$\llvert  J\rrvert  +\llvert   \partial(J)\rrvert   $ vertices. Assume
that such a forest is already
constructed. The events $\bigcap_{i \in T_l} \Omega _i$ are
independent for
different trees $T_l, l=1 ,\ldots,L$ in the forest. Hence,
\begin{eqnarray*}
\mathbb{P} \bigl(\llVert Wx \rrVert _2 \le t n^{-(c+1)/2}v t
\bigr) &\le& \mathbb{P} \bigl( \bigl\llvert (Wx)_i \bigr\rrvert \le
t n^{-(c+1)/2}v t \mbox{ for all } i \in[n] \bigr)
\\
&\le& \prod_{l=1}^L \mathbb{P} \biggl(
\bigcap_{i \in T_l} \Omega _i \biggr) \le
\prod_{l=1}^L t^{\llvert    T_l\rrvert   -1}
=t^{\llvert  J\rrvert  +\llvert
\partial(J)\rrvert   -L},
\end{eqnarray*}
where we used Lemma~\ref{l: tree} in the last inequality.
Since by the strong expansion condition, $\llvert  J\rrvert  +\llvert
\partial(J)\rrvert   -L \ge
(1+\kappa ) \llvert  J\rrvert  $, the last quantity is less than or equal to
$t^{(1+\kappa )
k}$ as required.

We proceed with the construction of the forest. At the first step, we
construct a spanning tree $\widetilde{T}_l$ for each connected component
of the set $J$. These trees are, obviously, disjoint, and $\sum_{l=1}^L
\llvert   \widetilde{T}_l\rrvert   =\llvert  J\rrvert  $. Now, we have to add
the vertices from
$\partial(J)$ as leaves to these trees. We do this by induction on $j
\in J$.
\begin{enumerate}[(2)]
\item[(1)] Let $j \in J$ be the smallest number. Add all vertices $i \in
\partial(J)$ connected to $j$ to the tree containing $j$ as the
descendants of $j$.

\item[(2)] Let $j \in J$ be the smallest number, which has not been used in
this process. Add all vertices $i \in\partial(J)$ connected to $j$,
which have not been already added, to the tree containing $j$ as its
descendants.
\end{enumerate}
Since any vertex in $\partial(J)$ is connected to some vertex in $J$,
the whole set $\partial(J)$ will be added at the end of this process.
Denote the trees obtained in this way by $T_1 ,\ldots,T_L$. The
construction guarantees that these trees are disjoint. This finishes
the construction of the forest and the proof of the lemma.
\end{pf}

Similarly to Lemma~\ref{l: very sparse}, we extend the small ball
probability result of Lemma~\ref{l: moderately compressible
individual} to a uniform bound using a net argument.
%
%le2.7 #&#

\begin{lemma} \label{l: moderately compressible uniform}
Let $W$ be an $n \times n$ skew-symmetric Gaussian matrix. Assume that
its large variances graph $\Gamma _W(n^{-c})$ satisfies the strong
expansion condition with parameter $\kappa \in(0,1)$ up to level $m<n$.
%MathOfer
Then there exists a constant $C'>0$ depending only on $c$ and $\kappa $
such that for any $t>0, k \le m$ and $v \in(0,1)$,
\begin{eqnarray*}
&& \mathbb{P} \bigl(\exists x \in\operatorname{Incomp}(k,v) \cap
\operatorname{Comp} \bigl((1+\kappa /2)k, \bigl(n^{-C'}v
\bigr)^{8/\kappa } \bigr):
\\
&&\qquad \llVert Wx \rrVert _2 \le n \cdot \bigl(n^{-C'}v
\bigr)^{8/\kappa } \mbox{ and }\llVert W \rrVert \le n \bigr) \le
e^{-k}.
\end{eqnarray*}
\end{lemma}

\begin{pf}
The proof repeats that of Lemma~\ref{l: very sparse}, so we only
sketch it.
For $t>0$, set
\[
\varepsilon=n^{-C-2}vt,
\]
where $C$ is the constant from Lemma~\ref{l: moderately compressible
individual}. Choose an $\varepsilon$-net $\mathcal{N}$ in
$\operatorname
{Sparse}((1+\kappa /2)k) \cap\operatorname{Incomp}(k,v)$ of cardinality
\[
\llvert \mathcal{N}\rrvert \le\pmatrix{n\cr (1+\kappa /2)k} \cdot \biggl(
\frac
{3}{\varepsilon} \biggr)^{(1+\kappa /2)k} \le \biggl( \frac{n^{\bar{c}}}{v t}
\biggr)^{(1+\kappa /2)k},
\]
where $\bar{c}$ depends only on $c$ and $\kappa $.
By the union bound,
\begin{eqnarray*}
\mathbb{P} \bigl(\exists x \in\mathcal{N}:\llVert Wx \rrVert _2 \le
n^{-C} v t\mbox{ and }\llVert W \rrVert \le n \bigr) &\le& \biggl(
\frac{n^{\bar{c}}}{v t} \biggr)^{(1+\kappa /2)k} \cdot t^{(1+\kappa )k}
\\
&\le& t^{(\kappa /4)k} \le e^{-k},
\end{eqnarray*}
provided that
\[
t = \biggl( \frac{v}{e n^{\bar{c}}} \biggr)^{4/\kappa }.
\]
Using an appropriately defined $C'>0$ depending only on $c$ and
$\kappa $,
and approximation by the points of the $\varepsilon$-net, we derive
from the previous inequality that
\begin{eqnarray*}
&& \mathbb{P} \bigl(\exists x \in\operatorname{Incomp}(k,v) \cap
\operatorname{Sparse} \bigl((1+\kappa /2)k \bigr):
\\
&&\qquad \llVert Wx \rrVert _2 \le2n \cdot \bigl(n^{-C'}v
\bigr)^{8/\kappa } \mbox{ and }\llVert W \rrVert \le n \bigr) \le
e^{-k}.
\end{eqnarray*}
To complete the proof, notice that for any vector $y \in\operatorname
{Incomp}(k,v) \cap\operatorname{Comp}((1+\kappa /2)k,
(n^{-C'}v)^{8/\kappa
})$, there is a vector $x \in\operatorname{Incomp}(k,v) \cap
\operatorname{Sparse}((1+\kappa /2)k) $ such that $\llVert  x-y
\rrVert _2 <
(n^{-C'}v)^{8/\kappa }$. The lemma now follows by again using approximation.
\end{pf}

Lemmas~\ref{l: very sparse} and~\ref{l: moderately compressible
uniform} can be combined to treat all compressible vectors. In the statement,
$d_0$ is a fixed, large enough universal positive integer.
%
%pr2.8 #&#

\begin{proposition} \label{p: all compressible}
Let $W$ be an $n \times n$ skew-symmetric Gaussian matrix.
Assume that its large variances graph $\Gamma _W(n^{-c})$ has minimal
degree $d\geq d_0$ and satisfies the strong expansion condition with parameter
$\kappa >0$ up to the level $m<n$.
%Set $\rho=(n/d)^{\phi(\k)}$.
Then there exists a constant $\phi(\kappa)$ depending on $\kappa$ only
so that, with $\rho=(n/d)^{\phi(\kappa )}$, one has
\[
\mathbb{P} \bigl(\exists x \in\operatorname{Comp} \bigl((1+\kappa /2)m,
n^{-\rho
} \bigr): \llVert Wx \rrVert _2 \le n^{-\rho+1}
\bigr) \le e^{-d/2}.
\]
\end{proposition}

%re2.9 #&#
%
\begin{remark} \label{r: phi(kappa)}
The proof below shows that it is enough to take
\[
\phi(\kappa )=\frac{c}{\kappa } \log \biggl( \frac{C}{\kappa } \biggr).
\]
\end{remark}

\begin{pf*}{Proof of Proposition \ref{p: all compressible}}
Set $v_0= n^{-c'}$, where $c'=\max(c,C')$, and $c, C'$ are the
constants from Lemmas~\ref{l: very sparse} and~\ref{l: moderately
compressible uniform}.
Let $L$ be the smallest natural number such that
\[
(d/2) (1+\kappa /2)^L \ge m.
\]
The definition of $L$ implies
\[
L \le\frac{\log(m/d)}{\log(1+\kappa /2)} \le\frac{c}{\kappa } \log\frac{n}{d}.
\]
For $l=1 ,\ldots,L-1$, define by induction
$v_{l+1}=(n^{-C'}v_l)^{8/\kappa }$,
where $C'$ is the constant from Lemma~\ref{l: moderately compressible uniform}.
The definition of $v_0$ implies that $v_l \le n^{-C'}$, so $v_{l+1} \ge
v_l^{16/\kappa }$, and thus,
\[
v_L \ge v_0^{(16/\kappa )^{L-1}} \ge n^{-\rho'}\qquad
\mbox{where } \rho'= \biggl( \frac{n}{d}
\biggr)^{(c/\kappa ) \cdot\log
(C/\kappa )}.
\]

We have
\begin{eqnarray*}
&& \operatorname{Comp} (m, v_L)
\\
&&\qquad \subset \operatorname{Comp} \bigl((d/2) (1+ \kappa /2)^L,
v_L \bigr)
\\
&&\qquad \subset \operatorname{Comp} (d/2, v_0)
\\
&&\quad\qquad{}\cup\bigcup_{l=1}^{L}
\operatorname{Comp} \bigl((d/2) (1+ \kappa /2)^l, v_l
\bigr) \setminus\operatorname{Comp} \bigl((d/2) (1+ \kappa /2)^{l-1},
v_{l-1} \bigr).
\end{eqnarray*}
Lemmas~\ref{l: very sparse} and~\ref{l: moderately compressible
uniform} combined with the union bound imply
\begin{eqnarray*}
&& \mathbb{P} \bigl(\exists x \in\operatorname{Comp} \bigl(m, n^{-\rho
'}
\bigr):
\llVert Wx \rrVert _2 \le n^{-\rho'+1}\mbox{ and }
\llVert W \rrVert \le n \bigr) \le e^{-d/2}.
\end{eqnarray*}
Applying Lemma~\ref{l: moderately compressible uniform} once more, we
derive the estimate
\begin{eqnarray*}
&& \mathbb{P} \bigl(\exists x \in\operatorname{Comp} \bigl((1+\kappa /2)m,
n^{-\rho
} \bigr):
\llVert Wx \rrVert _2 \le n^{-\rho+1}\mbox{ and }
\llVert W \rrVert \le n \bigr) \le e^{-d/2}
\end{eqnarray*}
with $\rho=(\rho')^{16/\kappa }$.

The proposition follows from the previous inequality and Lemma~\ref{l: norm}.
\end{pf*}

%s3 #&#
\section{The smallest singular value}
\label{sec: smallest singular}

The main result of this section is the following lower bound for the
smallest singular value of a Gaussian skew-symmetric matrix with a
strongly expanding large variances graph.
%
%th3.1 #&#

\begin{theorem} \label{th: smallest singular}
Let $n \in\mathbb{N}$ be an even number.
Let $V$ be an $n \times n$ skew-symmetric matrix, and denote by
$\Gamma $
its large variances graph $\Gamma _V(n^{-c})$.
Assume that:
\begin{longlist}[(3)]
\item[(1)] $\operatorname{Var}(v_{i,j}) \le1$ for all $i,j \in[n]$;
\item[(2)] the minimal degree of a vertex of $\Gamma $ is at least $d \ge2
\log n$;
\item[(3)]$\Gamma $ is $\kappa $-strongly expanding up to level $\frac
{n-d+2}{1+\kappa /4}$.
\end{longlist}
Then
\[
\mathbb{P} \bigl( s_n(V) \le t n^{-\tau } \bigr) \le n t
+e^{-cd},
\]
where $\tau =(n/d)^{\psi(\kappa )}$ for some positive $\psi(\kappa )$.
\end{theorem}

%re3.2 #&#
%
\begin{remark} \label{r: psi(kappa)}
Tracing the proof of Theorem~\ref{th: smallest singular} and using
Remark~\ref{r: phi(kappa)}, one can show that it is enough to take
\[
\psi(\kappa )= \frac{C'}{\kappa } \log \biggl( \frac{C}{\kappa } \biggr).
\]
\end{remark}

\begin{pf*}{Proof of Theorem~\ref{th: smallest singular}}
To prove the theorem, we use the negative second moment identity. Let
$A$ be an $n \times n$ matrix with columns $A_1 ,\ldots,A_n$. For $j
\in
[n]$, let $h_j \in S^{n-1}$ be a vector orthogonal to all columns of
$A$, except the $j$th one. Then
\[
\bigl\llVert A^{-1} \bigr\rrVert _{\mathrm{HS}}^2 = \sum
_{j=1}^n \bigl(h_j^T
A_j \bigr)^{-2}.
\]
Hence,
\[
s_n(A)=\frac{1}{\llVert  A^{-1} \rrVert } \ge\frac{1}{\llVert  A^{-1} \rrVert _{\mathrm{HS}}} \ge n^{-1/2}
\cdot\min_{j \in[n]} \bigl\llvert h_j^T
A_j \bigr\rrvert .
\]
Let $\rho$ be as in Proposition~\ref{p: all compressible}.
The argument above shows that if we use the matrix $V$ in place of $A$
and define the unit vectors $h_j, j \in[n]$ as before, then the
theorem would follow if the inequalities
%
%e3.1 #&#
%
\begin{equation}
\label{eq: h_j^T} \mathbb{P} \bigl( \bigl\llvert h_j^T
V_j \bigr\rrvert \le t n^{-\rho+c} \bigr) \le t
+e^{-cd}
\end{equation}
hold for all $j \in[n]$.
Indeed, the theorem follows from (\ref{eq: h_j^T}) and the assumption
on $d$ by the union bound.
We will establish inequality (\ref{eq: h_j^T}) for $j=1$. The other
cases are proved in the same way.

Let $W$ be the $(n-1) \times(n-1)$ block of $V$ consisting of rows and
columns from 2 to $n$. The matrix $W$ is skew-symmetric, and its large
variances graph is the subgraph of $\Gamma $ containing vertices $2
,\ldots,n$.
Therefore, $\Gamma $ has properties (1), (2), (3) with slightly relaxed
parameters. Indeed, property (1) remains unchanged. Property (2) is
valid with $d$ replaced by $d-1$. Property (3) is satisfied with
parameter $\kappa /2$ in place of $\kappa $ since for any $J \subset
\Gamma
\setminus\{1\}$, the boundary of $J$ in $\Gamma $ and in $\Gamma
\setminus\{
1\}$ differs by at most one vertex.

Recall that $W$ is a skew-symmetric matrix of an odd size. This matrix
is degenerate, so there exists $u \in S^{n-2}$ such that $Wu=0$.
This allows us to define the vector $h \in S^{n-1}$ orthogonal to the
columns $V_2 ,\ldots,V_n$ of the matrix $V$ by
\[
h= %
\pmatrix{0
\cr
u}.
\]
Define the event $\Omega $ by
\[
\Omega = \bigl\{ W: \exists u \in\operatorname{Comp} \bigl(n-d+1,
n^{-\rho} \bigr), Wu=0 \bigr\}.
\]
The graph $\Gamma _W(n^{-c}) \setminus\{1\}$ is $(\kappa /4)$ strongly
expanding up to level
$\frac{n-d+1}{1+\kappa /4}$.
By Proposition~\ref{p: all compressible}, $\mathbb{P}(\Omega ) \le e^{-d/2}$.
Condition on the matrix $W \in\Omega ^\complement$.
After the conditioning, we may assume that $u \in\operatorname
{Incomp}(n-d+1, n^{-\rho})$, and so the set $J = \{ j: \llvert
u_j\rrvert   \ge
n^{-\rho-1/2}\}$ has at least $n-d+1$ elements. Since the degree of
the vertex $\{1\}$ in the large variances graph of $V$ is at least $d$,
this means that there exists a $j \in J$ for which $\operatorname
{Var}(v_{j1})\ge n^{-c}$.
Therefore, conditionally on $W$, $h^TV_1=\sum_{j=2}^n u_j v_j$ is a
normal random variable with variance
\[
\operatorname{Var} \bigl(h^TV_1 \bigr) = \sum
_{j=2}^n u_j^2 \cdot
\operatorname{Var}(v_j) \ge n^{-2 \rho-1} \cdot n^{-c}.
\]
The bound on the density of a normal random variable implies
\[
\mathbb{P} \bigl( \bigl\llvert h^TV_1 \bigr\rrvert \le
C n^{-\rho-c/2-1/2} t | W \in\Omega ^\complement \bigr) \le t.
\]
Finally,
\begin{eqnarray*}
&&\mathbb{P} \bigl( \bigl\llvert h^TV_1 \bigr\rrvert
\le C n^{-\rho-c/2-1/2} t \bigr)
\\
&&\qquad\le\mathbb{P} \bigl( \bigl\llvert h^TV_1 \bigr
\rrvert \le C n^{-\rho-c/2-1/2} t \mid W \in\Omega ^\complement \bigr) +
\mathbb{P}( W \in \Omega )
\\
&&\qquad\le t +e^{-d/2}.
\end{eqnarray*}
This completes the proof of (\ref{eq: h_j^T}) for $j=1$. Since the
proof for the other values of $j$ is the same, it proves Theorem~\ref
{th: smallest singular}.
\end{pf*}
An immediate corollary of Theorem~\ref{th: smallest singular} is the following.
%
%co3.3 #&#

\begin{corollary}
\label{cor: matching exists}
Let $A$ be the adjacency matrix of a graph $\Gamma$ which satisfies:
\begin{longlist}[(2)]
\item[(1)] the minimal degree of a vertex of $\Gamma $ is at least $d>2\log n$;

\item[(2)] $\Gamma $ is $\kappa $-strongly expanding up to level $\frac
{n-d+2}{1+\kappa /4}$.
\end{longlist}
Then $A$ possesses a unique doubly stochastic scaling $B=DAD$ and the
graph $\Gamma$ possesses a perfect matching.
\end{corollary}

\begin{pf}{Proof} % of Corollary~\ref{cor: matching exists}}
We begin by showing that a perfect matching in~$\Gamma $ exists. Assume
otherwise. Then $s_n(A\odot G)=0$ since $\mathbb{E}\det(A\odot
G)=\operatorname{haf}
(A)=0$. The latter equality contradicts Theorem~\ref{th: smallest singular}.

To show that $A$ possesses a doubly stochastic scaling, choose an edge
$e=(u,v)$ in $\Gamma $ and create a graph $\Gamma'$ by erasing $u,v$, and
all edges attached to them from~$\Gamma$. The graph $\Gamma'$
satisfies assumptions (1) and (2) in the statement, with slightly
smaller constants $\kappa,d$. Thus, $\Gamma '$ possesses a perfect
matching. This implies that for any edge in $\Gamma $ there exists a
perfect matching containing that edge. By Bregman's theorem (\cite
{bregman}, Theorem~1), this implies that $A$ possesses a unique doubly
stochastic scaling $B=D_1AD_2$. The fact that $D_1=D_2$ follows from
the strict convexity of relative entropy and the characterization of
the doubly stochastic scaling as its minimizer; see \cite{bregman},
equation~(7).
\end{pf}

\section{Local bound on eigenvalues density}
\label{sec: local bounds}
In this section, we prove a general bound on the crowding of
eigenvalues at $0$ for a class of Hermitian matrices whose variance
matrix is doubly-stochastic. The results are somewhat more general than
our needs in the rest of the paper and may be of independent interest
and, therefore, we introduce new notation.

Let $X$ denote an $n\times n$ matrix, Hermitian (in\vspace*{1pt} the sense that
$X_{ij}^*=X_{ji}$), with entries $\{X_{ij}\}_{i\le j}$ that are
independent zero mean random variables. (In our application, the
$X_{ij}$ variables are all Gaussian.) Following
\cite{EKYY}, we set $s_{ij}=E\llvert    X_{ij}\rrvert   ^2$ and $\zeta
_{ij}=X_{ij}/\sqrt
{s_{ij}}$ (with $\zeta_{ij}=0$ if $s_{ij}=0$). We assume that
the variables $\zeta_{ij}$ possess uniformly bounded $p$ moments for
all $p>0$. Finally, we denote the eigenvalues of the matrix $X$ by
$\lambda_1(n)\geq\lambda_2(n)\geq\cdots\geq\lambda_n(n)$, and
use $L_n=n^{-1}\sum_{i=1}^n \delta_{\lambda_i(n)}$ for the empirical
measure of eigenvalues.

We assume that $\sum_j s_{ij}=1$, and, to avoid trivialities, that the
matrix $S=\{s_{ij}\}$ is irreducible (otherwise, the matrix $X$ can be
decomposed in blocks due to the symmetry).

Let $M=(\max_{ij} s_{ij})^{-1}$. We assume the following.
%
%as4.1 #&#

\begin{assumption}
\label{ass-M}
For some $\vartheta\in(0,1]$, one has that $M\geq n^{\vartheta}$.
\end{assumption}

With Assumption~\ref{ass-M}, we have the following proposition.
%
%pr4.2 #&#

\begin{proposition}
\label{prop-conc}
With notation and assumptions as in the setup above, fix $\epsilon
<1/5$. If Assumption~\ref{ass-M} holds then
for every $D>0$ there exists $n_0=n_0(\epsilon,D)$ such that for any
$n>n_0$, and with
\[
\mathcal{N}(\eta)= \bigl\llvert \bigl\{i: \lambda_i(n)\in(-\eta,
\eta) \bigr\} \bigr\rrvert ,
\]
one has
%
%e4.1 #&#
%
\begin{equation}
\label{eq-propconc} \mathbb{P} \bigl(\exists\eta\ge M^{-\epsilon}: \mathcal {N}(
\eta)> C n \cdot\eta \bigr)\leq n^{-D}.
\end{equation}
\end{proposition}

\begin{pf} We will use \cite{EKYY}, Theorem 2.3, a simplified form of
which we quote below after introducing some notation.
% We are interested in $E=0$ in their notation,
%and use the same notation $M$ as they do.
%Also following their notation,
Following the notation in \cite{EKYY}, we let $m(z)$ denote the
Stieltjes transform of the semicircle law, $m(z)=(-z+\sqrt{z^2-4})/2$,
and set
\[
\Gamma(z)= \bigl\llVert \bigl(1-m(z)^2S \bigr)^{-1} \bigr
\rrVert _{\ell^\infty\to\ell^\infty}.
\]
Note that with $z=i\eta$ one has by equation~(A.1) of \cite{EKYY}, that for some
universal constant $C$,
%
%e4.2 #&#
%
\begin{equation}
\label{eq-Gamma} \Gamma(i\eta)\leq\frac{C\log n}{\eta}.
\end{equation}
Introduce now, similarly to equation (2.14) of \cite{EKYY}
for a parameter
$\gamma>0$,
\begin{eqnarray*}
\tilde\eta &=& \min \biggl\{u>0: \frac{1}{Mu}\leq\min \biggl\{
\frac{M^{-\gamma}}{\Gamma(i\eta')^3}, \frac{M^{-2\gamma}}{\Gamma
(i\eta')^4 \cdot\Im(m(i\eta')) } \biggr\},
\\
&&{}  \mbox{for all }
\eta'\in[u,10] \biggr\}.
\end{eqnarray*}
(Note that we do not use $\widetilde\Gamma_N(z)$ as in \cite{EKYY} since
we only need the relation $\Gamma\geq\widetilde\Gamma_N$.)
Note that $\Im( m(i\eta'))=(\sqrt{4+\eta'^2}-\eta')/2$ is bounded
above and below by a universal constant for $\eta'\in[0,10]$. Hence,
using (\ref{eq-Gamma}), we get that
%
%e4.3 #&#
%
\begin{equation}
\label{eq-tildeeta} \tilde\eta\leq C \bigl(M^{2\gamma-1}(\log n)^4
\bigr)^{1/5}:=\bar\eta,
\end{equation}
for some universal constant $C$. For given $\varepsilon<1/5$, we will
chose $\gamma\in(0,1/2)$ and $n_0=n_0(\varepsilon, D)$ so that
$\bar\eta\le M^{-\varepsilon}$ whenever $n>n_0$.

Denote by $m_n(z)=n^{-1} \sum_{i=1}^n \frac{1}{\lambda_i-z}$
%\delta_{\lambda_i(n)}$
the Stieltjes transform of the empirical measure of eigenvalues of $X$.
We have the following.
%
%th4.3 #&#

\begin{theorem}[(\cite{EKYY}, Theorem 2.3)]\label{theo-EKYY}
For any $\gamma\in(0,1/2)$, any $\epsilon'>0$ and any $D>0$ there
exists an $n_0=n_0(D,\gamma,\epsilon')$ so that, uniformly in $\eta
\in[\tilde\eta,10]$, and for all $n>n_0$,
%
%e4.4 #&#
%
\begin{equation}
\label{eq-stieltjescontrol} \mathbb{P} \biggl( \bigl\llvert m_n(i\eta)- m(i\eta )
\bigr\rrvert > \frac{n^{\epsilon'}}{M\eta} \biggr)\leq n^{-D}.
\end{equation}
\end{theorem}

Fix $\eta\geq\bar\eta$.
Let $\mathcal{A}$ denote the complement of the event in (\ref
{eq-stieltjescontrol}).
Assume that $\mathcal{A}$ occurs.
Using the uniform boundedness of $m(i \eta)$ and inequality (\ref
{eq-tildeeta}), we obtain
\[
\Im m_n(i \eta)\leq C+\frac{n^{\epsilon'}}{M \eta}\leq C+C\frac
{n^{\epsilon'}}{M\cdot(M^{2\gamma-1}(\log n)^4)^{1/5}}.
\]
Choosing $\gamma$ and $\epsilon'$ small enough, we can guarantee that
the right-hand side in the last display is uniformly
bounded in $n$.
With such choice,
\[
C\geq\int\frac{ \eta}{\lambda^2+\eta^2} \,dL_n(\lambda) \geq\frac
{1}{2\eta}
\int_{- \eta}^{ \eta}\,dL_n(\lambda)=
\frac
{1}{2\eta}L_n \bigl([-\eta,\eta] \bigr)
\]
provided that $\mathcal{A}$ occurs.
This means that
\[
\mathbb{P} \bigl( \mathcal{N}(\eta)> C n\cdot\eta \bigr)\leq n^{-D}.
\]
To derive (\ref{eq-propconc}) from the previous inequality, one can
use the union bound over $\eta=1/k$ with $k \in\mathbb{N},
M^{\varepsilon
} \ge k \ge1$.\vadjust{\goodbreak}
\end{pf}

A better estimate can be obtained if one assumes a spectral gap. First,
we have the following.
% We begin with the following lemma.
%
%le4.4 #&#

\begin{lemma}
\label{lem-1draft}
$S$ has exactly one eigenvalue at $+1$ and at most one eigenvalue at $-1$.
\end{lemma}

\begin{pf} The claim concerning the eigenvalue at $1$ is the
Perron--Frobenius theorem.
To check the claim on the eigenvalues at $-1$, consider $S^2$. It may
be reducible, but at most to $2$ blocks. Indeed, suppose there are 3
disjoint blocks $A_1,A_2,A_3$, that is, disjoint subsets $A_i$ of
$[N]$, $i=1,2,3$, so that for all $a\in A_i, b\in A_j$ with
$i\neq j$ one has $S^2_{a,b}=0$. By the irreducibility of $S$, there is
a path of odd length connecting $A_1$ and $A_2$, and similarly there is
a path of odd length connecting $A_2$ and $A_3$. Hence, there is a path
of even length connecting $A_1$ and $A_3$, in contradiction with the
block disjointness of $S^2$. The claim now follows by applying the
Perron--Frobenius theorem to each of the blocks of $S^2$.
\end{pf}

%\begin{definition}
% \label{def-spectralgap}
% The matrix $S$ possesses a \textit{spectral gap} $\delta$ if there
% there do not exist eigenvalues of $S$ in $(-1,-1+\delta)\cup(1-
%\delta,1)$.
%\end{definition}
By Lemma~\ref{lem-1draft}, if $S$ has a spectral gap then
the eigenvalues at $1$ and $-1$ (if the later exists) are unique and isolated.
%
%pr4.5 #&#

\begin{proposition}
\label{prop-concgap}
With notation and assumptions as in the setup above, fix $\epsilon<1$.
If Assumption~\ref{ass-M} holds and $S$ possesses a spectral gap
$\delta$ then for every $D>0$ there exists $n_0=n_0(\epsilon,D,\delta
)$ such that for any $n>n_0$, with
\[
\mathcal{N}(\eta)= \bigl\llvert \bigl\{i: \lambda_i(n)\in(-\eta,
\eta) \bigr\} \bigr\rrvert ,
\]
one has
%
%e4.5 #&#
%
\begin{equation}
\label{eq-propconc1} \mathbb{P} \bigl(\exists\eta\ge M^{-\epsilon}: \mathcal {N}(
\eta)>n\cdot \eta \bigr)\leq n^{-D}.
\end{equation}
\end{proposition}

The proof is identical to that of Proposition~\ref{prop-conc}, using
\cite{AEK}, Theorem 1.1, instead of \cite{EKYY}, Theorem 1.2. We omit
the details.
%s5 #&#
\section{Concentration of the hafnian of a random matrix}
\label{sec: proof thm 1 2}
In this section, we prove Theorems~\ref{thm-1} and~\ref{thm-2}. Both
results follow from the concentration of the Gaussian measure for
Lipschitz functions. To this end, we consider a Gaussian vector
$G=(G_{ij})_{1 \le i<j \le n} \in\mathbb{R}^{n(n-1)/2}$ and use it to form
the skew-symmetric matrix $G^{\mathrm{skew}}$.
However, the function $F(G)=\log\operatorname{det}(\mathcal{B}\odot
G^{\mathrm{skew}})$
(where $\mathcal{B}_{ij}=\sqrt{B}_{ij}$)
is not Lipschitz. To overcome this obstacle, we write
\[
\log\operatorname{det} \bigl(\mathcal{B}\odot G^{\mathrm{skew}} \bigr)= \sum
_{j=1}^n \log s_j \bigl(
\mathcal{B}\odot G^{\mathrm{skew}} \bigr)
\]
and use Theorem~\ref{th: smallest singular} and Proposition~\ref
{prop-conc} to obtain lower bounds on the singular values
which are valid with probability close to 1. On this event,
we replace the function $\log$ by its truncated version, which makes
it Lipschitz with a controlled Lipschitz constant. Then an application
of the Gaussian concentration inequality yields the concentration of
the new truncated function about its expectation. This expectation is
close to
$\mathbb{E}\log\operatorname{det}(\mathcal{B}\odot G^{\mathrm{skew}})$.
Recall that instead of
the concentration about this value, we want to establish the concentration
about $\log\operatorname{haf}(B)=\log\mathbb{E}\operatorname
{det}(\mathcal{B}
\odot
G^{\mathrm{skew}})$. In other words, we have to
swap the expectation and the logarithm and estimate the error incurred
in this process. This will be achieved due to the fast decay of the
tail in the concentration inequality.

\begin{pf*}{Proof of Theorem~\ref{thm-1}}
The proof proceeds as in \cite{RZ}, Section~7. (The argument can be
traced back to \cite{FRZ}.)
Without loss of generality, we may assume that $n>n_0$, where
$n_0=n_0(\varepsilon,D)$ appears in Proposition~\ref{prop-conc}.
Indeed, if $n \le n_0$, we can
choose the constant $C$ in the formulation of the theorem appropriately
large, so that $C n_0^{1-\varepsilon\vartheta} \ge n_0$.
In this case, Theorem~\ref{thm-1} follows from Barvinok's theorem.

Fix $\epsilon<1/5$, $D>4$ as in the statement of the theorem, and
$t=n^{-(D+1)}$. With $\tau=\tau(\kappa,\alpha )$ as in Theorem~\ref{th:
smallest singular}
(with $X$ replaced by $W$)
and $\mathcal{N}(\eta)=\mathcal{N}(\eta)(W)$ as in Proposition~\ref
{prop-conc},
introduce the events
\[
\mathcal{W}_1= \bigl\{s_n(W)\leq tn^{-\tau}
\bigr\}, \qquad\mathcal{W}_2= \bigl\{\mathcal{N} \bigl(n^{-\varepsilon
\vartheta}
\bigr)\geq n^{1-\epsilon\vartheta} \bigr\},\qquad\mathcal{W}=\mathcal
{W}_1\cup\mathcal{W}_2.
\]
By Theorem~\ref{th: smallest singular} and Proposition~\ref
{prop-conc}, we have that for all $n>n_0$,
%
%e5.1 #&#
%
\begin{equation}
\label{eq-conc0} \mathbb{P}(\mathcal{W})\leq3n^{-D}.
\end{equation}
Let $\widetilde\det(W)=\prod_{i} (\llvert   \lambda_i(W)\rrvert   \vee
n^{-\epsilon
\vartheta})$. Note that on $\mathcal{W}^\complement$ we have that
%
%e5.2 #&#
%
\begin{equation}
\label{eq-conc3} \bigl\llvert \log\widetilde\det(W)-\log\det (W) \bigr\rrvert \leq
C(D)n^{1-\epsilon
\vartheta} \log n.
\end{equation}

Set $U=\log\widetilde\det(W)-\mathbb{E}\log\widetilde\det(W)$.
We next derive concentration results for $U$.
The map $(\lambda_i(W))_{i=1}^N
\to\log\widetilde\det(W)$ is Lipschitz with constant
$n^{1/2+\epsilon\vartheta}$. Therefore,
by standard concentration for
the Gaussian distribution (see \cite{GZ,Ledoux}),
using that the variance of the entries of $W$ is bounded above
by $n^{-\vartheta}$,
we have for some universal constant $C$ and any $u>0$,
%
%e5.3 #&#
%
\begin{eqnarray}
\label{eq-conc1} \mathbb{P} \bigl(\llvert U\rrvert >u \bigr) &=& \mathbb {P}
\bigl( \bigl\llvert \log \widetilde \det(W)-\mathbb{E}\log\widetilde\det(W) \bigr
\rrvert >u \bigr)
\nonumber\\[-8pt]\\[-8pt]\nonumber
&\leq& \exp \biggl(- \frac{Cu^2}{n^{1+(2\epsilon-1)\vartheta}} \biggr).
\nonumber
\end{eqnarray}
Therefore,
%
%e5.4 #&#
%
\begin{equation}
\label{eq-070714a} \mathbb{E} \bigl(e^{\llvert    U\rrvert   } \bigr)\leq 1+\int
_0^\infty \exp \biggl(u- \frac
{Cu^2}{n^{1+(2\epsilon-1)\vartheta}}
\biggr) \,du \leq\exp \bigl( n^{1+(2\epsilon-1)\vartheta} \bigr).
\end{equation}
In particular, we obtain that
%
%e5.5 #&#
%
\begin{equation}
\label{eq-conc070714} \mathbb{E}\log\widetilde\det{W}\leq\log \mathbb{E} \widetilde\det
{W}\leq\mathbb{E}\log \widetilde\det{W} +n^{1+(2\epsilon-1)\vartheta}.
\end{equation}
The first inequality above follows from Jensen's inequality, and the
second one from (\ref{eq-070714a}).

We can now complete the proof of the theorem.
We have by Markov's inequality that
%
%e5.6 #&#
%
\begin{equation}
\label{eq-070714} \mathbb{P} \bigl(\log\det(W)-\log\mathbb{E}\det
(W)>n^{1-\epsilon
\vartheta}\log n \bigr) \leq e^{-n^{1-\epsilon\vartheta}\log n}.
\end{equation}
On the other hand, note that $\mathbb{E}\det(W)\leq\mathbb
{E}\widetilde\det(W)$.
Therefore,
with $C(D)$ as in~(\ref{eq-conc3}),
\begin{eqnarray*}
&& \mathbb{P} \bigl(\log\det(W)-\log\mathbb{E}\det(W)\leq- \bigl(C(D)+2
\bigr)n^{1-\epsilon
\vartheta
}\log n \bigr)
\\
&&\qquad\leq\mathbb{P} \bigl(\log\det(W)-\log\mathbb {E}\widetilde\det(W)\leq -
\bigl(C(D)+2 \bigr)n^{1-\epsilon\vartheta}\log n \bigr)
\\
&&\qquad\leq\mathbb{P} \bigl(\log\widetilde\det(W)-\log\mathbb {E}\widetilde
\det(W) \leq-2n^{1-\epsilon\vartheta}\log n \bigr)+P \bigl(\mathcal {W}^\complement
\bigr),
\end{eqnarray*}
where (\ref{eq-conc3}) was used in the last display.
Using now (\ref{eq-conc0}) and
the upper bound in~(\ref{eq-conc070714}), we get
\begin{eqnarray*}
\label{eq-070714b} &&\mathbb{P} \bigl(\log\det(W)-\log\mathbb {E}\det(W)\leq -
\bigl(C(D)+2 \bigr)n^{1-\epsilon
\vartheta}\log n \bigr)
\\
&&\qquad\leq3n^{-D}+ \mathbb{P} \bigl(\log\widetilde\det (W)-
\mathbb{E}\log \widetilde \det(W)\leq-2n^{1-\epsilon\vartheta}\log n+n^{1+(2\epsilon
-1)\vartheta}
\bigr).
\end{eqnarray*}
Using that $\epsilon<1/5$ and applying
(\ref{eq-conc1}), we conclude that
\[
\mathbb{P} \bigl(\log\det(W)-\log\mathbb{E}\det(W)\leq- \bigl(C(D)+2
\bigr)n^{1-\epsilon
\vartheta
}\log n \bigr) \leq4n^{-D}.
\]
Together with (\ref{eq-070714}), it yields
\[
\mathbb{P} \bigl( \bigl\llvert \log\det(W)-\log\mathbb{E}\det (W) \bigr\rrvert
\leq \bigl(C(D)+2 \bigr)n^{1-\epsilon
\vartheta}\log n \bigr) \leq5 n^{-D}.
\]
To obtain the statement of the theorem, we prove the previous
inequality with $\varepsilon' \in(\varepsilon,1/5)$ instead of
$\varepsilon$, and then choose $C'>0$ such that $(C(D)+2)
n^{1-\varepsilon\vartheta} \log n \le C' n^{1-\varepsilon'
\vartheta}$.
%Together with \eqref{eq-070714}, this completes the proof of the
%theorem,
%once $\epsilon$ is adjusted in order to get rid of the $\log n$ factor.
%
%On the other hand, writing $A:=\left\left\vert  \log\widetilde\det(W)-\log
%\det(W)\right\left\vert  $ we
%have
%\begin{eqnarray}
%\label{eq-conc4}
%\left\left\vert  \E\log\widetilde\det(W)-\E\log\det(W)\right\left\vert
%&\leq&\E A=
%\E\left(A{\mathbf1}_{\WW^\complement}\right)+
%\E\left(A{\mathbf1}_{\WW}\right)\nonumber\\
%&\leq&C(D)n^{1-\epsilon\vartheta}\log n+C(3n^{-D})^{1/2} n^{3/2},
%\end{eqnarray}
%where in the last display we used \eqref{eq-conc1}, the Cauchy-Schwarz
%inequality,
%\eqref{eq-conc0} and Lemma~\ref{lem-7.2}. In particular, we get
%that for any $n>n_0(\epsilon,D,c,\vartheta)$,
%\begin{equation}
%\label{eq-conc6}
%\left\right\vert  \E\log\widetilde\det(W)-\E\log\det(W)\right\right\vert  \leq\E A
%\leq2 C(D)n^{1-
%\epsilon\vartheta}\log n.
%\end{equation}
%Together with \eqref{eq-conc1} and \eqref{eq-conc0}, this leads to
%$$\P(\left\right\vert  \log\det(W)-\E\log\det(W)\right\right\vert  >C'(D)n^{1-
%\epsilon\vartheta}\log n)
%\leq3n^{-D}.\quad$$
%Modifying $\epsilon$ if needed one can get rid of the $\log n$ factor.
\end{pf*}
The proof of Theorem~\ref{thm-2} is similar to that of Theorem~\ref{thm-1}. However,
to exploit the tighter bounds on the intermediate singular values
provided by Proposition~\ref{prop-concgap}, we use a different
truncation, redefining the function $\log\widetilde\det(\cdot)$ and
estimate its Lipschitz constant more accurately.

\begin{pf*}{Proof of Theorem~\ref{thm-2}}
Fix $\varepsilon\in(1/2,1)$. As in the proof of Theorem~\ref
{thm-1}, we may assume that $n>n_0(\varepsilon,D,\delta)$, where
$n_0(\varepsilon,D,\delta)$ was introduced in Proposition~\ref
{prop-concgap}.
The inequality (\ref{eq-propconc1}) can be rewritten as
%
%e5.7 #&#
%
\begin{equation}
\label{eq-sing values} \mathbb{P} \biggl(\exists k \ge n^{1-\varepsilon
\vartheta
}:
s_{n-k}(W) \le c \frac{k}{n} \biggr) \le n^{-D}.
\end{equation}
This inequality can be used to bound the Lipschitz constant of the
truncated logarithm. Let
%
%e5.8 #&#
%
\begin{equation}
\label{eq: m_0 large} m_0 \ge n^{1-\varepsilon\vartheta}
\end{equation}
be a number to be chosen later. Define
%
%e5.9 #&#
%
\begin{equation}
\label{eq: log tilde det} \log\widehat{\operatorname{det}}(W) = \sum
_{k=1}^n \phi_k \bigl(s_{n-k}(W)
\bigr) = \sum_{k=1}^n \log
\bigl(s_{n-k}(W) \vee\varepsilon_k \bigr),
\end{equation}
where
\[
\varepsilon_k= \cases{ \displaystyle c \frac{m_0}{n}, &\quad for $k
<m_0$,
\vspace*{3pt}\cr
\displaystyle c \frac{k}{n}, &\quad for $k \ge
m_0$.}
\]
Denote for a moment $N=n(n-1)/2$. For a vector $Y \in\mathbb{R}^N$ consider
the $n \times n$ skew symmetric matrix $Y^{\mathrm{skew}}$ whose entries above
the main diagonal equal to the corresponding entries of $Y$. Let
$\mathcal{B}$
be the $n \times n$ matrix whose entries are square roots of the
corresponding entries of $B$.
Note that the function $F: \mathbb{R}^N \to\mathbb{R}$ defined by
\[
F(G)=\log\widehat{\operatorname{det}} \bigl(\mathcal{B}\odot G^{\mathrm{skew}}
\bigr)
\]
is the composition of three functions: $F =F_3 \circ F_2 \circ F_1$, where:
\begin{longlist}[(2)]
\item[(1)] $F_1: \mathbb{R}^N \to\mathbb{R}^{n^2}, F_1(G)=\mathcal
{B}\odot G^{\mathrm{skew}}$, whose
Lipschitz constant does not exceed $n^{-\vartheta/2}$;

\item[(2)] $F: \mathbb{R}^{n^2} \to\mathbb{R}^n_+$ defined by $F(W)=
(s_1(W) ,\ldots,
s_n(W) )$, which is \mbox{1-}\break Lipschitz;

\item[(3)]$F_3: \mathbb{R}^n_+ \to\mathbb{R}, F_3(x_1 ,\ldots,x_n)=\sum_{k=1}^n \phi
_k(x_{n-k})$, where $\phi_k$ is defined in~(\ref{eq: log tilde det}).
\end{longlist}
By the Cauchy--Schwarz inequality,
\begin{eqnarray*}
\llVert F_3 \rrVert _{\mathrm{Lip}} &\le& \Biggl( \sum
_{k=1}^n \llVert \phi_k \rrVert
_{\mathrm{Lip}}^2 \Biggr)^{1/2} \le \biggl( \biggl(
\frac{n}{c m_0} \biggr)^2 \cdot m_0 + \sum
_{k>m_0} \biggl( \frac{n}{c k} \biggr)^2
\biggr)^{1/2}
\\
&\le& C \frac{n}{\sqrt{m_0}}.
\end{eqnarray*}
Therefore,
\[
\llVert F \rrVert _{\mathrm{Lip}} \le C \frac{n^{1-\vartheta
/2}}{\sqrt{m_0}}.
\]
Applying the standard Gaussian concentration for Lipschitz functions,
we obtain
%
%e5.10 #&#
%
\begin{eqnarray}
\label{eq-tild-conc} \mathbb{P} \bigl( \bigl\llvert \log\widehat {\operatorname
{det}}(W)- \mathbb{E} \log\widehat{\operatorname{det}}(W) \bigr\rrvert \ge u
\bigr) &\le& 2 \exp \biggl(- \frac{u^2}{2 \llVert  F \rrVert _{\mathrm{Lip}}^2} \biggr)
\nonumber\\[-8pt]\\[-8pt]\nonumber
&\le& 2 \exp \bigl(- c m_0 n^{\vartheta-2} u^2
\bigr)
\nonumber
\end{eqnarray}
which replaces formula (5.8) in the proof of Theorem~\ref{thm-1}.
Arguing as in the proof of that theorem, we obtain
%
%e5.11 #&#
%
\begin{equation}
\label{eq-wt} \mathbb{E}\log\widehat\det(W)\leq\log\mathbb {E}\widehat\det (W)
\leq\mathbb{E}\log \widehat\det(W) + C_0 \frac{n^{2-\vartheta}}{m_0}
\end{equation}
from the inequality above.

Let $C'=D+\tau (\kappa ,\alpha )$, where $\tau (\kappa ,\alpha )$
is as in Theorem~\ref
{th: smallest singular}.
Set
\[
\mathcal{W}_1= \bigl\{s_n(W)\leq n^{-C'}
\bigr\}, \qquad\mathcal{W}_2= \biggl\{ \exists k \ge n^{1-\varepsilon
\vartheta},
s_{n-k}(W) \le c \frac{k}{n} \biggr\},
\]
and let $\mathcal{W}=\mathcal{W}_1\cup\mathcal{W}_2$.
Then Theorem~\ref{th: smallest singular} and (\ref{eq-sing values})
imply $\mathbb{P}(\mathcal{W}) \le n^{-D}$. On $\mathcal
{W}^\complement$ we have
%
%e5.12 #&#
%
\begin{equation}
\label{eq-tilde-no-tilde} \bigl\llvert \log\widehat{\operatorname {det}}(W)- \log
\operatorname{det}(W) \bigr\rrvert \le C m_0 \log n,
\end{equation}
which plays the role of (5.7).
Arguing as in the proof of Theorem~\ref{thm-1}, we show that
\begin{eqnarray*}
&&\mathbb{P} \biggl( \log\det(W)- \log\mathbb{E}\det(W) \le- \biggl(C
m_0 \log n+ 2C_0 \frac{n^{2-\vartheta}}{m_0} \biggr) \biggr)
\\
&&\qquad\le\mathbb{P} \biggl( \log\widehat\det(W)- \log\mathbb {E}\widehat\det
(W) \le- 2C_0 \frac{n^{2-\vartheta}}{m_0} \biggr) +\mathbb{P}(\mathcal{W})
\\
&&\qquad\le\mathbb{P} \biggl( \log\widehat\det(W)- \mathbb {E}\log\widehat\det
(W) \le- C_0 \frac{n^{2-\vartheta}}{m_0} \biggr) +\mathbb{P}(\mathcal{W})
\\
&&\qquad\le\exp \biggl( -C' \frac{n^{2-\vartheta}}{m_0} \biggr) +
n^{-D} \le2 n^{-D}.
\end{eqnarray*}
Here, the first inequality follows from (\ref{eq-tilde-no-tilde}) and
$\log\mathbb{E}\det(W) \le\log\mathbb{E}\widehat\det(W)$, the
second one from
the upper bound in (\ref{eq-wt}), and the third one from (\ref{eq-tild-conc}).

We select the optimal $m_0=C n^{1-\vartheta/2} \log^{-1/2} n$ in the
inequality above. Since $\varepsilon>1/2$, condition (\ref{eq: m_0
large}) holds for sufficiently large $n$.
\[
\mathbb{P} \bigl( \log\det(W)- \log\mathbb{E}\det(W) \le- C n^{1-\vartheta/2}
\log^{1/2} n \bigr) \le2 n^{-D}.
\]
Combining this with the bound
\[
\mathbb{P} \bigl( \log\det(W)- \log\mathbb{E}\det(W) \ge C n^{1-\vartheta/2}
\log^{1/2} n \bigr) \le n^{-D}
\]
following from Markov's inequality, we complete the proof.
\end{pf*}

%s6 #&#
\section{Doubly stochastic scaling and proof of Theorem \texorpdfstring{\protect\ref
{thm-3bis}}{1.5}}
\label{sec: doubly stochastic}
To prove Theorem~\ref{thm-3bis}, we have to scale the adjacency matrix
of the graph in order to apply Theorems~\ref{thm-1} and~\ref{thm-2}.
The existence of such scaling has been already established in Corollary
\ref{cor: matching exists}. We will show now that the smallest
nonzero entry of the scaled adjacency matrix is at least polynomial in
$n$. This crucial step in the proof of Theorem~\ref{thm-3bis} allows
us to conclude that the large entries graph of the scaled matrix
coincides with the original graph.
%
%pr6.1 #&#

\begin{proposition}
\label{prop: doubly stochastic}
Fix $\alpha ,\kappa>0$. Let $A$ be the adjacency matrix of a graph
$\Gamma$ whose minimal degree satisfies $d\geq\alpha n+2$. Assume that
%
%e6.1 #&#
%
\begin{equation}
\label{eq-SEpropds} \mbox{$\Gamma$ is $\kappa$ strongly expanding up to level $n(1-
\alpha )/(1+\kappa/4)$.}
\end{equation}
Then there exists a constant $\nu$ so that $A$ possesses a doubly
stochastic scaling $B=DAD$ with\vspace*{-3pt}
\[
\min_i D_{ii}\geq n^{-\nu}.\vspace*{-3pt}
\]
\end{proposition}

In particular, under the assumptions of Proposition~\ref{prop: doubly
stochastic}, we have that
%
%e6.2 #&#
%
\begin{equation}
\label{eq-minimal D} B_{ij}\geq n^{-2\nu}\qquad\mbox{whenever }
B_{ij}>0.
\end{equation}

Before describing the proof of Proposition~\ref{prop: doubly
stochastic}, let us state a complementary claim, which says that under
stronger expansion conditions on $\Gamma$ we can guarantee that the
entries in its scaled adjacency matrix are \textit{polynomially small}.
This ensures that
%Ofer #5
under this stronger expansion property,
condition (\ref{eq-dsadj}) in Theorem~\ref{thm-3bis} is automatically
satisfied.

%Ofer #6
In what follows, if $X$ is a set of vertices in a graph then
$E(X,X)$ denotes those edges in the graph connecting vertices in $X$.
%
%pr6.2 #&#

\begin{proposition}
\label{prop:not-large}
Fix $\alpha, \kappa, \varepsilon> 0$. There exists a constant
$\theta$ depending
only on $\alpha$ and $\kappa$ such that the following holds.

Let $A$ be the adjacency matrix of a
graph $\Gamma$ whose minimal degree satisfies $d\geq\alpha n+2$.
%Ofer #7
Assume that for any subset $X$ of vertices satisfying
$\llvert    X\rrvert   \ge\frac{d}{4}$ and $\llvert    E(X,X) \rrvert   \le
\theta\cdot
n^{1+\varepsilon}$, it
holds
that
\[
\llvert \partial_s X\rrvert \ge(1 + \kappa) \cdot\llvert X\rrvert
,
\]
where $\partial_s X$ denotes the set of external neighbors of $X$ such
that any $y \in\partial_s X$ has at least $\frac{d \llvert    X\rrvert   }{10n}$
neighbors in $X$.

Then $A$ possesses a doubly stochastic scaling $B=DAD$ with
\[
\max_i D_{ii}\leq n^{-\varepsilon/2}.
\]
\end{proposition}

In particular, under the assumptions of Proposition~\ref
{prop:not-large}, we have that
%
%e6.3 #&#
%
\begin{equation}
\label{eq-minimal Dnew} B_{ij}\leq n^{-\varepsilon}\qquad\mbox{whenever }
B_{ij}>0.
\end{equation}

To prove Proposition~\ref{prop: doubly stochastic}, we argue by
contradiction. Assume that one of the diagonal entries, say $D_{11}$,
is smaller than $n^{-\nu}$, where $\nu=\nu(\alpha ,\kappa )$ will
be chosen
at the end of the proof. The double stochasticity of the scaled matrix
implies that there exists a neighbor $i \sim1$ for which the
corresponding entry of the scaling matrix
$D$ is large. In fact, we can prove this for more than one entry. In
Lemma~\ref{lem:big-big}, we construct a set $X=X_0$ of vertices of
cardinality at least $d/2$ such that the corresponding entries of the
scaling matrix are greater than $(1/2)n^{\nu-2}$. We use this as a
base of induction. In Lemma~\ref{lem:bigger}, we show that there
exists a set $X_1$ of vertices of cardinality $\llvert    X_1\rrvert   \ge
(1+\beta)
\llvert    X_0\rrvert   $ containing $X_1$ such that all entries of\vadjust{\goodbreak} the
scaling matrix
corresponding to $X_1$ are still polynomially large. Proceeding by
induction, we construct an increasing sequence of sets $X_0 \subset X_1
\subset\cdots\subset X_l$ such that $\llvert    X_l\rrvert   \ge(1+\beta
)^l \llvert    X_0\rrvert   $,
and all diagonal entries corresponding to the vertices of $X_l$ are
greater than 1. The number of induction steps $l$ which we are able to
perform will depend on $\nu$. If $\nu$ is chosen large enough, then
we will get $(1+\beta)^l \llvert    X_0\rrvert   >n$, reaching the desired
contradiction.

The proof of Proposition~\ref{prop:not-large} is very similar. Assume,
toward contradiction, that, say, $D_{nn}$ is larger than
$n^{-\varepsilon/2}$. By the double stochasticity of the scaled
matrix, there exists a set
$\mathcal{A}$ of neighbors $i \sim n$, for which the corresponding
entries of
the scaling matrix $D$ are small.
%Ofer #8
Using again the double stochasticity of the scaled matrix
produces a set $X=X_0$ of vertices of cardinality at least $d/4$ such
that the corresponding entries of the scaling matrix are greater than
$(\alpha/8) n^{-\varepsilon/2}$. We use this as induction base. In
Lemma~\ref{lem:step}, we show that there exists a set $X_1$ of
vertices of cardinality $|   X_1|  \ge(1+\gamma) |
X_0 |  $ containing $X_1$
such that all entries of the scaling matrix corresponding to $X_1$ are
still large. Proceeding by induction, we construct an increasing
sequence of sets $X_0 \subset X_1 \subset\cdots\subset X_l$ such that
$\llvert    X_l\rrvert   \ge(1+\gamma)^l \llvert    X_0\rrvert   $, and all
diagonal entries corresponding
to the vertices of $X_l$ are greater than $\Omega
(n^{-\varepsilon/2} )$. The number of induction steps $l$ which
we are able to perform will depend on $\theta$. If $\theta= \theta
(\alpha,\kappa)$ is chosen large enough, then we will get $(1+\gamma
)^l \llvert    X_0\rrvert   >n$, reaching contradiction.

\begin{pf*}{Proof of Proposition~\ref{prop: doubly stochastic}}
Without loss of generality, we assume throughout that the constants
$\alpha
$ and $\kappa$ are small enough so that
%
%e6.4 #&#
%
\begin{equation}
\label{eq-smallc} (1-\alpha )/(1-\kappa/4)> 1/2.
\end{equation}

By Corollary~\ref{cor: matching exists}, $A$ possesses a doubly
stochastic scaling $B=DAD$, where $D=\operatorname{Diag}(r_1,\ldots
,r_n)$. Without loss of generality, we assume that $r_1\leq r_2\leq
\cdots\leq r_n$.
Note that since $B$ is doubly stochastic,
%
%e6.5 #&#
%
\begin{equation}
\label{eq-req} r_i= \biggl(\sum_{j\sim i}
r_j \biggr)^{-1}.
\end{equation}
We will need a few simple lemmas.
%Ofer #9
%
%le6.3 #&#

\begin{lemma}
\label{lem:Markov}
Let $s_1,\ldots,s_n \in[0,1]$ and assume that
$\sum_{i=1}^n s_i \ge S$. Then, for any $0 < \gamma< 1$, there exists a
subset $I \subseteq[n]$ of cardinality at least $(1-\gamma) \cdot S$,
%such that for each $i \in I$ holds
$s_i \ge\gamma\cdot S/n$ for each $i\in I$.
\end{lemma}

\begin{pf} Assume otherwise. Then there are at least $n-(1-\gamma)S$
elements $s_i<\gamma\cdot(S/n)$. Therefore, %
\[
\sum_{i=1}^n s_i\leq(1-
\gamma)S+ \bigl(n-(1-\gamma)S \bigr)\cdot\gamma\frac{S}{n}<S.
\]\upqed
\end{pf}

The next lemma quantifies the following intuition: given a large set
$\mathcal{A}$ of indices corresponding to small
entries of the scaling matrix, we can find a large set of indices
(neighbors of $\mathcal{A}$) corresponding to large entries of the
scaling matrix.
%
%le6.4 #&#

\begin{lemma}
\label{lem:big-from-small}
%Ofer #9A Fixed notation screwup
Let $\mathcal{A}\subseteq[n]$ such that
$r_i\leq\mu$
for all $i \in\mathcal{A}$.
Then, for any $0 < \gamma< 1$, there exists a subset $X \subseteq[n]$,
of cardinality at least $(1-\gamma) \cdot\llvert   \mathcal{A}\rrvert   $, such
that for all
$j \in X$,
\[
r_j \ge\frac{\gamma}{\mu n}.
\]
\end{lemma}

\begin{pf}
Denote by $B = (b_{ij} )$ the doubly stochastic scaling of $A$.
For $1 \le i \le n$, let $s_i = \sum_{j \in\mathcal{A}} b_{ij}$.
By the double-stochasticity of $B$, we have $\sum_{i=1}^n s_i =
\sum_{j \in\mathcal{A}} \sum_{i=1}^n b_{ij} = \llvert   \mathcal
{A}\rrvert   $.

By Lemma~\ref{lem:Markov}, there is a set $X$ of indices with
$\llvert    X\rrvert   \ge(1-\gamma) \cdot\llvert   \mathcal{A}\rrvert   $ with $s_i
\ge\gamma\cdot
\frac{\llvert   \mathcal{A}\rrvert   }{n}$ for each $i \in X$. Let $i \in X$.
We have
\[
r_i = \frac{s_i}{\sum_{j \in\mathcal{A}, j \sim i} r_j} \ge\gamma \cdot \frac
{\llvert   \mathcal{A}\rrvert   /n}{\llvert   \mathcal{A}\rrvert   \cdot\mu} =
\frac{\gamma}{\mu n}.
\]\upqed
\end{pf}

The next lemma is the base of our inductive construction.
%
%le6.5 #&#

\begin{lemma}
\label{lem:big-big}
Let $Q = 1/r_1$. Then there exists a subset $X$ of $[n]$ of cardinality
at least $d/2$ such that for each $i \in X$,
\[
r_i \ge\frac{Q}{2n^2}.
\]
\end{lemma}

\begin{pf}
By (\ref{eq-req}), $\sum_{i \sim1} r_i = 1/r_1 = Q$. Therefore,
there is at least one index $i_0 \sim1$ for which $r_{i_0} \ge Q/n$.
Let $\mathcal{A}$ be the set of neighbors of $i_0$. Then $ \llvert
\mathcal{A}\rrvert
\geq d$ and
for all $j \in\mathcal{A}$, $r_j \le n/Q$.

The proof is completed by an application of Lemma~\ref
{lem:big-from-small} with $\gamma= 1/2$.
\end{pf}

Lemma~\ref{lem:bigger} below will be used for the inductive step.
%
%le6.6 #&#

\begin{lemma}
\label{lem:bigger}
Let $m \ge n$. Let $X$ be a subset of indices such that $r_i\geq m$ for
each $i \in X$. Then:
\begin{longlist}[(2)]
\item[(1)]
$\llvert   \partial X\rrvert   \ge(1 + \kappa) \cdot\llvert    X\rrvert   $;

\item[(2)]
there exists a subset $Z$ of indices, disjoint from $X$, of cardinality
at least $(\kappa\cdot\llvert    X\rrvert   -1)/2$ such that each $j \in
Z$ satisfies
\[
r_j \ge\frac{\kappa\llvert    X\rrvert   -1}{2n^2} \cdot m.
\]
\end{longlist}
\end{lemma}

\begin{pf}
Clearly, no two vertices in $X$ are connected in $\Gamma $ (otherwise, we
would have an entry of size at least $m^2 > 1$ after scaling).
Therefore, $X$ is a set of disconnected vertices, and, since $\Gamma $
contains a perfect matching, we have $\llvert    X\rrvert   \le n/2$.

%Ofer #10
Since $X$ is disconnected,
%This, by
(\ref{eq-SEpropds}) and (\ref{eq-smallc}) imply that
%, means that $\left\left\vert   X\right\right\vert  \le\frac{n-d}{1+\kappa/4}$, and
%therefore, since
%$X$ is disconnected,
%
\[
\llvert \partial X\rrvert \ge(1 + \kappa) \cdot\llvert X\rrvert ,
\]
proving the first claim of the lemma.

Let $Y:= \partial X$. We note that $r_i\leq1/m$ for any $i \in Y$.
To show the second claim of the lemma, we will find a subset $Z$ of
indices, disjoint with $X \cup Y$, such that $r_i \ge\frac{\kappa
\llvert    X\rrvert   -1}{2n^{2}} \cdot m$ for all $i \in Z$.

Let $C = (X \cup Y)^\complement$. Recall that $B$ is the doubly
stochastic scaling of $A$. Since
\[
\sum_{i\in C\cup X\cup Y} \sum_{j\in Y}
b_{ij}= \sum_i \sum
_{j\in Y} b_{ij}=\llvert Y\rrvert ,
\]
we have
\begin{eqnarray*}
\sum_{i \in C, j \in Y} b_{ij} &=& \llvert Y\rrvert -
\sum_{i \in X, j \in Y} b_{ij} - \sum
_{i \in Y, j \in Y} b_{ij}
\\
%Ofer #11 I only see inequality!
& \geq& \bigl(\llvert Y\rrvert - \llvert X\rrvert \bigr) - \sum
_{i \in Y, j \in Y} b_{ij}
\\
&\ge& \bigl(\llvert Y\rrvert -
\llvert X\rrvert \bigr) - \frac{n^2}{m^2}
\\
&\ge& \kappa\cdot\llvert X\rrvert - 1,
\end{eqnarray*}
where in the second
inequality we used that $r_i\leq1/m$ for $i\in Y$.

For $i \in C$, set $s_i = \sum_{j \in Y} b_{ij}$. Then $\sum_{i \in
C} s_i \ge\kappa\cdot\llvert    X\rrvert   -1$. By Lemma~\ref{lem:Markov},
there is a
set of at least $(\kappa\cdot\llvert    X\rrvert   -1)/2$ indices $i$, for which
\[
s_i \ge\frac{\kappa\llvert    X\rrvert   -1}{2 \llvert    Y\rrvert   } \ge \frac{\kappa\llvert    X\rrvert   -1}{2n}.
\]
Call this set $Z$.
For each $i \in Z$, we have
\[
r_i = \frac{s_i}{\sum_{j \in Y, j\sim i} r_j} \ge\frac{(\kappa
\llvert    X\rrvert   -1)/(2n)}{\llvert    Y\rrvert   /m} \ge
\frac{\kappa\llvert    X\rrvert   -1}{2n^2} \cdot m,
\]
completing the proof of the lemma.
\end{pf}
%
%Ofer #12
We are
now ready to perform the inductive procedure
proving Proposition~\ref{prop: doubly stochastic}. Let
\[
R = \biggl(\frac{4}{\kappa\alpha} \cdot n \biggr)^{c/(\kappa
\log(1/\alpha))},
\]
for a sufficiently large $c$. We will assume that $r_1 < 1/R$, and
reach a contradiction.

We use Lemma~\ref{lem:big-big} to construct a set $X$ of cardinality
at least $d/2$ such that $r_i \ge R/(2n^2)$.
for all $i \in X$.

Assuming $m:= R/(2n^2) \ge n$, which we may, we can now apply
Lemma~\ref{lem:bigger} to construct a set $Z$ disjoint from $X$, of
cardinality at least $(\kappa\cdot\llvert    X\rrvert   -1)/2$ such that
\[
r_j \ge\frac{\kappa\llvert    X\rrvert   -1}{2n^2} \cdot m \ge\frac{\kappa
\alpha }{4n} \cdot m
\qquad\mbox{for all }j \in Z.
\]

We now define $X_0:= X$, $m_0:= m$, $Z_0:= Z$; and set $X_1 = X_0
\cup Z_0$, $m_1 = \frac{\kappa\alpha }{4n} \cdot m_0$, and apply
Lemma~\ref{lem:bigger} to $X_1$ (assuming $m_1$ is not too small). We
continue this process to obtain an increasing sequence of sets
$X_0,X_2,\ldots,X_t$. Since
\[
n \ge\llvert X_t\rrvert \ge(1 + \kappa/2) \cdot\llvert
X_{t-1} \rrvert \ge\cdots \ge(1 + \kappa/2)^t \cdot
X_0 \ge(1 + \kappa/2)^t \cdot\frac{\alpha}{2} \cdot
n,
\]
the number of steps $t$ is upper bounded by $c_1 \cdot\kappa\log
1/\alpha$, for some absolute constant~$c_1$. On the other hand, if $c$
in the definition of $R$ is large enough, the number of steps will be
larger than that, reaching a contradiction.
\end{pf*}

\begin{pf*}{Proof of Proposition~\ref{prop:not-large}}
Assume, for contradiction's sake, that $r_n \ge n^{-\varepsilon/2}$.
Since $\sum_{i \sim n } r_i = 1/r_n \le n^{\varepsilon/2}$,
for at least half of neighbors of $n$ holds $r_i \le2\cdot
n^{\varepsilon/2 - 1}$.
%Ofer #13
Let $\mathcal{A}$ be the set of these neighbors. By our
assumption on the minimal degree in $\Gamma$, we have $\llvert   \mathcal{A}
\rrvert   \ge
\frac{d}{2}$.

Applying Lemma~\ref{lem:big-from-small} to $\mathcal{A}$ with $\gamma
= 1/2$
and $\mu= 2\cdot n^{\varepsilon/2 - 1}$ gives a subset $X_0$ of
$[n]$ of cardinality at least $\frac{d}{4}$ such that for all $i \in
X_0$ holds $r_i \ge\alpha/8 \cdot n^{-\varepsilon/2}$. This is our
induction base.

An inductive step is provided by the following lemma.
%
%le6.7 #&#

\begin{lemma}
\label{lem:step}
Fix a constant $b \ge1/\sqrt{\theta}$. Let $X \subseteq[n]$ be such
that $\llvert    X\rrvert   \ge\frac{\alpha}{4} \cdot n$ and for any $i
\in X$, $r_i
\ge b \cdot n^{-\varepsilon/2}$. Then\vspace*{1pt} there exists a subset $X'$ of
$[n]$ of cardinality at least $ (1 + \frac{\kappa(1-\kappa
)}{2} ) \cdot\llvert    X\rrvert   $ such that for all $j \in X'$ it holds that
$r_j \ge\frac{\alpha^2 \kappa b}{80} \cdot n^{-\varepsilon/2}$.
\end{lemma}

\begin{pf}
Since $n \ge\sum_{i,j \in X} d_{ij} \ge\llvert    E(X,X)\rrvert   \cdot b^2
n^{-\varepsilon}$, we have
\[
\bigl\llvert E(X,X) \bigr\rrvert \le b^{-2} \cdot n^{1+\varepsilon}
\le \theta\cdot n^{1+\varepsilon}.
\]
Hence, by our assumptions on the graph $\Gamma$, we have $\llvert
\partial_s
X\rrvert   \ge(1 + \kappa) \cdot\llvert    X\rrvert   $. For each $j \in
\partial_s X$ holds
\[
r_j \le\frac{1}{\sum_{i \in X,i \sim j} r_i} \le\frac{1}{
(\alpha/10 ) \cdot\llvert    X\rrvert   \cdot b n^{-\varepsilon/2}} \le
\frac
{40}{b \alpha^2} \cdot n^{\varepsilon/2 - 1}.
\]

Applying Lemma~\ref{lem:big-from-small} with $\gamma= \kappa/2$ and
$\mu= \frac{40}{b \alpha^2} \cdot n^{\varepsilon/2 - 1}$ to
$\mathcal{A}=
\partial_s X$, produces a set $X'$ satisfying the requirements of the lemma.
\end{pf}

We now ready to perform the inductive procedure proving Proposition
\ref{prop:not-large}. Let
\[
R = \biggl(\frac{1}{\kappa} \cdot\log \biggl(\frac{4}{\alpha
} \biggr)
\biggr)^{c/(\alpha^2 \kappa)},
\]
for a sufficiently large $c$. We will assume that $\theta> 1/R$, and
reach a contradiction.

We start constructing the sequence $\{X_i\}$, starting from the set
$X_0$ constructed above, and applying Lemma~\ref{lem:step}
iteratively. Clearly, we should stop after at most $S = \log_{1 +
(\kappa(1-\kappa)/2)} (\frac{4}{\alpha} )$
steps. However, if $c$ in the definition of $R$ is large enough, we
would be able to make more steps than that, reaching a contradiction.
\end{pf*}

We now combine the bound (\ref{eq-smallc}) on the scaled matrix with
Theorems~\ref{thm-1} and~\ref{thm-2} to derive Theorem~\ref{thm-3bis}.

\begin{pf*}{Proof of Theorem~\ref{thm-3bis}}
Recall that $B=DAD$ denotes the doubly stochastic scaling of $A$ and
$G^{\mathrm{skew}}$ denotes a skew symmetric matrix with independent $N(0,1)$
entries above the main diagonal. Note that
\[
\det \bigl(A\odot G^{\mathrm{skew}} \bigr)=\frac{1}{\det(D)} \det
\bigl(B_{1/2}\odot G^{\mathrm{skew}} \bigr),
\]
where $B_{1/2}$ denotes the matrix whose entries are the square roots
of the entries of $B$. Therefore, it is enough to consider the
concentration for $\det(B_{1/2}\odot G^{\mathrm{skew}})$.
The proof of Theorem~\ref{thm-3bis} now follows by applying Theorems
\ref{thm-1} and~\ref{thm-2}.
\end{pf*}
%
%s7 #&#
\section{The strong expansion condition}
\label{sec: strong expansion}
As noted in the \hyperref[sec1]{Introduction},
the strong expansion condition is stronger than the classical vertex
expansion condition
\[
\forall J \subset[M] \qquad\llvert J\rrvert \le M/2\quad\Rightarrow\quad \bigl\llvert
\partial(J) \bigr\rrvert \ge\kappa \llvert J\rrvert .
\]
It might have been desirable to replace the strong expansion property
by a weaker and more natural classical vertex expansion condition.
Proposition
\ref{p: bias} from the \hyperref[sec1]{Introduction} shows that not only the latter
condition is insufficient to guarantee a subexponential error in
Barvinok's estimator, but in fact there is an example of a graph $G$
with associated random matrix $W$ that barely misses the strong
expansion property, for which Barvinok's estimator yields an
exponential error with high probability. We provide here the proof of
Proposition~\ref{p: bias}.
%\begin{proposition} \label{p: bias}
% Let $\d>0$. For any $N \in\N$, there exists a graph $\G$ with $M>N$
%vertices %such that
% \begin{equation} \label{c: weak expansion}
% \forall J \subset[M] \quad\left\left\vert   J\right\right\vert  \le M/2 \Rightarrow
%% \left\left\vert  \partial(J)\right\left\vert  - (1-\d) \left\left\vert  \mbox{Con}(J)
%\right\right\vert  \ge\k\left\right\vert   J\right\right\vert
% \end{equation}
% and
% \[
% \P\left(\frac{\det(W_{\G})}{\E\det(W_{\G})} \le e^{-cM} \right)
% \ge1- e^{-c'M}.
%\]
% Here $c,c', \k$ are constants depending on $\d$.
%\end{proposition}
%
\begin{pf*}{Proof of Proposition~\ref{p: bias}}
Without loss of generality, assume that $\delta<1/6$.
Let $n \in\mathbb{N}$. Set $m=\lfloor\frac{\delta n }{2}\rfloor$.
Define a graph $\Gamma $ with $M=2(m+n)$ vertices as follows.
\begin{itemize}
\item The vertices in $[n]$ form a clique, which will be called the center.
\item Any of the vertices in $[n+1: 2(n+m)]$, called peripheral, is
connected to all vertices of the center.
\item In addition, for $k>n$, the vertices $2k-1$ and $2k$ are
connected to each other. (See Figure~\ref{fig1}.)
\end{itemize}
%
%%%%%%%%%%%%%%%%%%%%%%%%%%%%%%%%%%%%%%%%%%%%%%%%%%%%%%%%%%%%%%%%%%%%%%%%

The adjacency matrix of $\Gamma $ has the block shape
\[
\pmatrix{ \mathbf{Q}_{n \times n} & & \mathbf{1}_{n \times n} & & \mathbf
{1}_{n \times2m}
\cr
\mathbf{1}_{n \times n} & & \mathbf{0}_{n \times n}
& & \mathbf{0}_{n \times2m}
\cr
& & & \bolds{\Delta } & &
\cr
\mathbf{1}_{2m \times n} & & \mathbf{0}_{2m \times n} & & \ddots
\cr
& & & & &
\bolds{\Delta }}.
\]
Here, $\mathbf{Q}_{n \times n}$ is the adjacency matrix of the
$n$-clique, that is, the matrix with 0 on the main diagonal and 1
everywhere else; $\mathbf{1}_{k \times l}$ is the $k \times l$ matrix
whose entries are equal to 1, and $\bolds{\Delta} $ is a $2 \times2$ matrix
\[
\bolds{\Delta} = \pmatrix{ 0 & 1
\cr
1 & 0}.
\]
The right lower block of $\Gamma $ contains $m$ such matrices $\Delta
$ on the
main diagonal.

\begin{figure}

\includegraphics{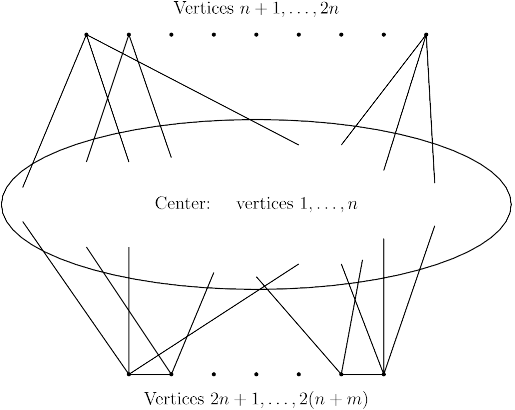}

\caption{Graph $\Gamma$.}\label{fig1}
\end{figure}

The matrix $W_{\Gamma }$ has the similar form
\[
W_{\Gamma }= \pmatrix{ \widetilde{\bolds{Q}}_{n \times n} & &
\mathbf{G}_{n \times n} & & \mathbf{G'}_{n \times2m}
\cr
-
\mathbf{G}_{n \times n} & & \mathbf{0}_{n \times n} & & \mathbf
{0}_{n \times2m}
\cr
& & & \widetilde{\bolds{\Delta }_1} & &
\cr
- \mathbf{G'}_{2m \times n} & & \mathbf{0}_{2m \times n} & &
\ddots
\cr
& & & & & \widetilde{\bolds{\Delta }}_m},
\]
where $\widetilde{\bolds{Q}}_{n \times n}$ is the $n \times n$
skew-symmetric Gaussian matrix; $\mathbf{G}_{n \times n}, \mathbf
{G'}_{n \times2m}$ are independent Gaussian matrices, and
\[
\widetilde{\bolds{\Delta }}_k= \pmatrix{ 0 & g_k
\cr
-g_k & 0}
\]
with independent $N(0,1)$ random variables $g_1 ,\ldots,g_m$.
Recall that
\[
\mathbb{E}\det W_{\Gamma } = \# \operatorname{Matchings}(\Gamma )
\]
the number of perfect matchings of the graph $\Gamma $.
Any vertex from $[n+1:2n]$ has to be matched to a vertex from the
center, which can be done in $n!$ ways.
Hence, for $k>n$, any vertex $2k-1$ has to be matched to its peripheral
neighbor $2k$, which can be done in the unique way.
Thus,
\[
\# \operatorname{Matchings}(\Gamma ) =n!>0.
\]
Consider $ \det W_{\Gamma }$.
Let $c>0$ be a constant to be chosen later.
A simple pigeonhole argument shows that
\[
\operatorname{det} W_{\Gamma }= F \bigl(\mathbf{G}_{n \times n},
\mathbf {G'}_{n \times2m} \bigr) \cdot\prod
_{j=1}^m g_j^2,
\]
where $F(\mathbf{G}_{n \times n}, \mathbf{G'}_{n \times2m})$ is a
homogeneous\vspace*{1pt} polynomial of degree $2n$ of entries of $\mathbf{G}_{n
\times n}$ and $\mathbf{G'}_{n \times2m}$.
Hence, for $\alpha = \frac{4}{\delta}+1$, we have
\begin{eqnarray*}
\mathbb{P} \biggl(\frac{\operatorname{det}(W_{\Gamma })}{\mathbb
{E}\operatorname
{det}(W_{\Gamma })} \ge e^{-cM} \biggr) &=&\mathbb{P}
\biggl( \frac{\operatorname{det}(W_{\Gamma })}{\mathbb{E}\operatorname
{det}(W_{\Gamma })} \ge\exp(- c \alpha m ) \biggr)
\\
&\le& \mathbb{P} \biggl( \frac{F(\mathbf{G}_{n \times n}, \mathbf{G'}_{n
\times2m})}{ \mathbb{E}F(\mathbf{G}_{n \times n}, \mathbf{G'}_{n
\times
2m})} \ge\exp( c \alpha m ) \biggr)
\\
&&{}+ \mathbb{P} \biggl(\frac{\prod_{j=1}^m g_j^2}{\mathbb{E}\prod_{j=1}^m g_j^2} \ge \exp(- 2c \alpha m ) \biggr).
\end{eqnarray*}
The first term above is smaller than $\exp(- c \alpha m )$ by
the Chebyshev inequality.
The second term also does not exceed $\exp(-c'm)$ if the constant $c$
is chosen small enough.

This proves the part of the proposition related to the error of the
Barvinok estimator.

It remains to check that condition (\ref{c: weak expansion}) is satisfied.
Let $J \subset[M]$ be a set of cardinality $\llvert  J\rrvert  \le
M/2$. If $J$
contains a vertex from the center, then $\llvert   \operatorname
{Con}(J)\rrvert   =1$
and $\partial(J)= [M] \setminus J$, so condition (\ref{c: weak
expansion}) holds.

Assume that $J \cap[n]= \varnothing$. Then $\llvert   \operatorname
{Con}(J)\rrvert
\le m+n = M/2$. Also, $\partial(J) \supset[n]$, so
\[
\bigl\llvert \partial(J) \bigr\rrvert \ge n \ge\frac{1}{1+\delta/2} \cdot
\frac{M}{2}.
\]
Therefore, since $\delta< 1/6$,
\[
\bigl\llvert \partial(J) \bigr\rrvert - (1-\delta) \bigl\llvert
\operatorname{Con}(J) \bigr\rrvert \ge \biggl( \frac{1}{1+\delta/2} -(1-\delta)
\biggr) \cdot\frac{M}{2} \ge \frac{\delta}{8} \cdot\frac{M}{2} \ge
\kappa \cdot\llvert J\rrvert
\]
if we choose $\kappa =\delta/8$.
This completes the proof of the proposition.
\end{pf*}

%\begin{appendix}
%\section{}
%\end{appendix}

\section*{Acknowledgment}
We thank Alexander Barvinok
for many helpful discussions.%s2 #&#

%\begin{supplement}[id=suppA]
%\slink[doi]{10.1214/00-AOPXXXXSUPP} %[doi,text={...}] - jei reikia
%suskaldyti doi
%\sfilename{aopXXXX\_supp.pdf}
%\end{supplement}

% imsref loaded by linak, 2015-08-11 14:12:45
%

\printaddresses
\end{document}